\documentclass[11pt]{article}
\usepackage{amsfonts}
\usepackage{graphics}
\usepackage{indentfirst}
\usepackage[colorlinks, linkcolor=red]{hyperref}
\usepackage{color, xcolor}
\hypersetup{urlcolor=red, citecolor=red}
\usepackage{cite}
\usepackage{latexsym}
\usepackage[paper=a4paper, left=1.6cm, right=1.6cm, top=1.8cm, bottom=1.6cm, headheight=5.5pt, footskip=0.8cm, footnotesep=0.8cm, centering, includefoot]{geometry}
\usepackage{amsmath}
\allowdisplaybreaks
\usepackage{amssymb}
\usepackage[dvips]{epsfig}
\usepackage{amscd}

\newtheorem{theorem}{Theorem}[section]
\newtheorem{remark}{Remark}[section]

\newtheorem{lemma}{Lemma}[section]

\DeclareMathOperator{\divv}{div}

\makeatletter
\@addtoreset{equation}{section}
\makeatother
\makeatletter
\@addtoreset{equation}{section}
\makeatother

\title{Global well-posedness and exponential decay of 2D nonhomogeneous Navier-Stokes and magnetohydrodynamic equations with density-dependent viscosity and vacuum
\thanks{This research was partially supported by National Natural Science Foundation of China (Nos. 11901474, 12071359), Exceptional Young Talents Project of Chongqing Talent (No. cstc2021ycjh-bgzxm0153), and the Innovation Support Program for Chongqing Overseas Returnees (No. cx2020082).}
}

\author{ Xin Zhong\thanks{School of Mathematics and Statistics, Southwest University, Chongqing 400715,
People's Republic of China ({\tt xzhong1014@amss.ac.cn}).
}
}
\date{ }

\begin{document}
\maketitle

\begin{abstract}
We establish global well-posedness of strong solutions for the nonhomogeneous magnetohydrodynamic equations with density-dependent viscosity and initial density allowing vanish in two-dimensional (2D) bounded domains. Applying delicate energy estimates and Desjardins' interpolation inequality, we derive the global existence of a unique strong solution provided that $\|\nabla\mu(\rho_0)\|_{L^q}$ is suitably small. Moreover, we also obtain exponential decay rates of the solution.
In particular, there is no need to impose some compatibility condition on the initial data despite the presence of vacuum.
As a direct application, it is shown that the similar result also holds for the nonhomogeneous Navier-Stokes equations with density-dependent viscosity.
\end{abstract}



\section{Introduction}
Magnetohydrodynamics studies the dynamics of electrically conducting fluids
and the theory of the macroscopic interaction of electrically conducting fluids with a magnetic field. The dynamic motion of the fluid and the magnetic field interact strongly with each other, so the hydrodynamic and electrodynamic effects are coupled. In the present paper,
let $\Omega\subset\mathbb{R}^2$ be a bounded smooth domain, we are concerned with nonhomogeneous magnetohydrodynamic equations with density-dependent viscosity in $\Omega$:
\begin{align}\label{1.1}
\begin{cases}
\rho_{t}+\divv(\rho\mathbf{u})=0,\\
(\rho\mathbf{u})_{t}+\divv(\rho\mathbf{u}\otimes\mathbf{u})
-\divv(2\mu(\rho)\mathfrak{D}(\mathbf{u}))+\nabla P=\mathbf{b}\cdot\nabla\mathbf{b}, \\
\mathbf{b}_t-\nu\Delta\mathbf{b}+\mathbf{u}\cdot\nabla\mathbf{b}
-\mathbf{b}\cdot\nabla\mathbf{u}=\mathbf{0},\\
\divv\mathbf{u}= \divv\mathbf{b}=0,
\end{cases}
\end{align}
with the initial condition
\begin{equation}\label{1.2}
(\rho,\rho\mathbf{u},\mathbf{b})(x,0)=(\rho_0,\rho_0\mathbf{u}_0,\mathbf{b}_0)(x),\ \ x\in\Omega,
\end{equation}
and the Dirichlet boundary condition
\begin{equation}\label{1.3}
(\mathbf{u},\mathbf{b})=(\mathbf{0},\mathbf{0}),\ x\in\partial\Omega,\ t>0.
\end{equation}
Here $\rho, \mathbf{u}, \mathbf{b}, P$ are the fluid density, velocity, magnetic field, and pressure, respectively. $\mathfrak{D}(\mathbf{u})$ denotes the deformation tensor given by
\begin{equation*}
\mathfrak{D}(\mathbf{u})=\frac{1}{2}(\nabla\mathbf{u}+(\nabla\mathbf{u})^{tr}).
\end{equation*}
The viscosity coefficient $\mu(\rho)$ is a function of the density satisfying \begin{equation}\label{1.4}
\mu\in C^1[0,\infty),\ \ \mu\geq\underline{\mu}>0
\end{equation}
for some positive constant $\underline{\mu}$, while the constant $\nu>0$ is the magnetic diffusion coefficient.

Notice that if there is no electromagnetic field effect (i.e., $\mathbf{b}=\mathbf{0}$), \eqref{1.1} reduces to the nonhomogeneous Navier-Stokes equations with variable viscosity
\begin{align}\label{1.6}
\begin{cases}
\rho_{t}+\divv(\rho\mathbf{u})=0,\\
(\rho\mathbf{u})_{t}+\divv(\rho\mathbf{u}\otimes\mathbf{u})
-\divv(2\mu(\rho)\mathfrak{D}(\mathbf{u}))+\nabla P=\mathbf{0}, \\
\divv\mathbf{u}=0.
\end{cases}
\end{align}
Many authors dealt with the above system. Lions \cite[Chapter 2]{L1996} derived the global weak solutions, yet the uniqueness and regularities of such weak solutions are big open questions. Later,
Desjardins \cite{D1997} introduced the so-called \textit{pesudo-energy method} and established global weak solutions with higher regularity for 2D
case provided that $\|\mu(\rho_0)-1\|_{L^\infty}$ is suitably small.
It should be noted that the solution obtained by Desjardins \cite{D1997} still does not have uniqueness. The main difficulty is due to the fact that Riesz transform does not map continuously
from $L^\infty$ to $L^\infty$ (see \cite{D1997} for details).
Meanwhile, if the initial density belongs to some
Besov spaces with positive index which guarantee that the initial density is at least a continuous function, Abidi and Zhang \cite{AZ2015} can show the uniqueness of the solution in the whole plane. Recently, some attention was focused on the well-posedness of strong solutions to \eqref{1.6}. For the initial density strictly away from vacuum, Abidi and Zhang \cite{AZ20152} proved the global well-posedness to the 3D Cauchy problem of \eqref{1.6} under the smallness assumptions on both $\|\mathbf{u}_0\|_{L^2}\|\nabla\mathbf{u}_0\|_{L^2}$ and $\|\mu(\rho_0)-1\|_{L^\infty}$. On the other hand, when the initial density allows vacuum, under the compatibility condition
\begin{equation}\label{CK}
-\divv(2\mu(\rho_0)\mathfrak{D}(\mathbf{u}_0))+\nabla P_0=\sqrt{\rho_0}\mathbf{g}\ \ \text{for some}\ (P_0,\mathbf{g})\in H^1\times L^2,
\end{equation}
which is proposed by Cho and Kim \cite{CK2004} in order to obtain the local existence of solutions solutions (see also \cite{LS2019} for an improved result), Huang and Wang \cite{HW2015} and Zhang \cite{Z2015} proved the global existence of strong solutions of \eqref{1.6} in 3D bounded domains provided the initial velocity is suitably small in some sense. Huang and Wang \cite{HW2014} also obtained the global strong solutions in 2D bounded
domains. Very recently, by time weighted techniques and energy methods, He et al. \cite{HLL2021} and Liu \cite{L2021} established global well-posedness of strong solutions to the 3D Cauchy problem without using the compatibility condition \eqref{CK} under suitable smallness conditions. Moreover, they also obtained exponential decay rates of the solution.

Let's turn our attention to the study of nonhomogeneous magnetohydrodynamic equations with variable viscosity coefficient. On one hand, in the absence of vacuum, Abidi and Paicu \cite{AP2008} obtained the global wellposedness
of strong solutions to the 3D Cauchy problem in the critical Besov space under
the assumptions that the initial velocity and magnetic field are small enough and the initial density $\rho_0$ approaches a positive constant. Sokrani \cite{S2020} proved global existence of strong solutions when the initial
data are small in some Sobolev spaces (see also related work \cite{SY2016}), which generalized the result for
nonhomogeneous Navier-Stokes equations with variable viscosity obtained by Abidi and Zhang \cite{AZ20152}. On the other hand, for the initial density allowing vacuum states, under the compatibility condition
\begin{equation}\label{LHY}
-\divv(2\mu(\rho_0)\mathfrak{D}(\mathbf{u}_0))+\nabla P_0-\mathbf{b}_0\cdot\nabla\mathbf{b}_0=\sqrt{\rho_0}\mathbf{g}\ \ \text{for some}\ (P_0,\mathbf{g})\in H^1\times L^2,
\end{equation}
Li \cite{L2018} showed global-in-time unique strong solution with density-dependent viscosity and resistivity coefficients to the 3D case under the condition that $\|\nabla\mathbf{u}_0\|_{L^2}+\|\nabla\mathbf{b}_0\|_{L^2}$ is small enough. This result was later improved by Liu \cite{L2019} without using \eqref{LHY} provided that $\|\rho_0\|_{L^\infty}+\|\mathbf{b}_0\|_{L^3}$
is suitably small (see related work \cite{L20212} for 2D case). Very recently, Zhang \cite{Z2020} investigated the global existence and large
time asymptotic behavior of strong solutions to the 3D Cauchy problem provided that the initial velocity and
magnetic field are suitable small in the $\dot{H}^\beta$-norm for some $\beta\in(\frac12,1]$.

As pointed out in \cite{HLL2021,Z2015}, the strong interaction between density and velocity will bring some serious difficulties in the mathematical study of global theory for the case of density-dependent viscosity, and the methods used for the case of constant viscosity cannot be applied directly.
In the present paper, we aim at investigating the global existence and exponential decay of strong solutions to the problem \eqref{1.1}--\eqref{1.3} without using the compatibility condition \eqref{LHY} via time weighted techniques.

Before stating our main result, we first explain the notations and conventions used throughout this paper. For $1\leq p\leq\infty$ and integer $k>0$, we use $L^p=L^p(\Omega)$ and $W^{k,p}=W^{k,p}(\Omega)$ to denote the standard Lebesgue and Sobolev spaces, respectively. When $p=2$, we use $H^k=W^{k,2}(\Omega)$. The space $H_{0,\sigma}^{1}$ stands for the closure in $H^1$ of the space $C_{0,\sigma}^\infty:=\{\pmb{\phi}\in C_{0}^\infty(\Omega)|\divv\pmb{\phi}=0\}$.

Our main result reads as follows:
\begin{theorem}\label{thm1.1}
Let the initial data $(\rho_0\geq0,\mathbf{u}_0,\mathbf{b}_0)$ satisfy
\begin{align}\label{A}
\rho_0\in W^{1,q}(\Omega),\ (\mathbf{u}_0,\mathbf{b}_0)\in H_{0,\sigma}^{1}(\Omega)\times H_{0,\sigma}^{1}(\Omega),\ q\in(2,\infty).
\end{align}
Then there exists a small positive constant $\varepsilon_0$ depending only on $\Omega,\underline{\mu},\bar{\mu}:=\sup\limits_{[0,\|\rho_0\|_{L^\infty}]}
\mu(\rho), \nu, q, \|\rho_0\|_{L^\infty}$, $\|\nabla\mathbf{u}_0\|_{L^2}^2$, and $\|\nabla\mathbf{b}_0\|_{L^2}^2$ such that if
\begin{align}\label{C}
\|\nabla\mu(\rho_0)\|_{L^q}\leq\varepsilon_0,
\end{align}
the problem \eqref{1.1}--\eqref{1.3} has a unique global strong solution
$(\rho\geq0,\mathbf{u},\mathbf{b})$ satisfying for $\tau>0$ and $r\in(2,q)$,
\begin{align}\label{1.5}
\begin{cases}
\rho\in L^\infty(0,\infty;W^{1,q})\cap C([0,\infty);W^{1,q}), \\
\mathbf{u}\in L^\infty(0,\infty;H^1)\cap L^\infty(\tau,\infty;H^{2})\cap L^2(\tau,\infty;W^{2,r}),\\
\nabla P\in L^\infty(\tau,\infty;L^2)\cap L^2(\tau,\infty;L^r),\\
\mathbf{b}\in L^\infty(0,\infty;H^1)\cap L^\infty(\tau,\infty;H^{2})\cap L^2(\tau,\infty;H^3),\\
\nabla\mathbf{u},\nabla\mathbf{b}\in C([\tau,\infty);L^2),\
\rho\mathbf{u}, \mathbf{b}\in C([0,\infty);L^2),\\
t\sqrt{\rho}\mathbf{u}_t,\ t\mathbf{b}_t
\in L^\infty(0,\infty;L^2),\\
e^{\frac{\sigma}{2}t}\nabla\mathbf{u},\
e^{\frac{\sigma}{2}t}\nabla\mathbf{b},\ e^{\frac{\sigma}{2}t}\sqrt{\rho}\mathbf{u}_t,\
e^{\frac{\sigma}{2}t}\Delta\mathbf{b}\in L^2(0,\infty;L^2),\\
t\nabla\mathbf{u}_t,\ t\nabla\mathbf{b}_t\in L^2(0,\infty;L^2),
\end{cases}
\end{align}
where $\sigma:=\min\left\{\frac{\underline{\mu}}{d^2\|\rho_0\|_{L^\infty}},
\frac{\nu}{d^2}\right\}$ with $d$ being the diameter of $\Omega$.
Moreover, there exists a positive constant $C$ depends only on $\Omega,\underline{\mu},\bar{\mu}, \nu, q, \|\rho_0\|_{L^\infty}, \|\nabla\mathbf{u}_0\|_{L^2}^2$, and $\|\nabla\mathbf{b}_0\|_{L^2}^2$ such that for $t\geq1$,
\begin{align}\label{lv1.2}
\|\mathbf{u}(\cdot,t)\|_{H^2}^2
+\|\nabla P(\cdot,t)\|^2_{L^2}+\|\mathbf{b}(\cdot,t)\|_{H^2}^2
+\|\sqrt{\rho}\mathbf{u}_t\|_{L^2}^2+\|\mathbf{b}_t\|_{L^2}^2
\le Ce^{-\sigma t}.
\end{align}
\end{theorem}


\begin{remark}
The conclusion in Theorem \ref{thm1.1} is somewhat surprising since the smallness condition \eqref{C} is independent of the initial magnetic field explicitly and just the same as that of nonhomogeneous Navier-Stokes equations (see \cite{HW2014}), which is in sharp contrast to the recent works \cite{L2019,L2018,L20212,Z2020}, where the authors considered global strong solution to the nonhomogeneous magnetohydrodynamic equations with density-dependent viscosity involving small initial magnetic field in some sense.
\end{remark}

\begin{remark}
Since the viscosity $\mu(\rho)$ depends on $\rho$, in order to bound the $L^2$-norm of the gradients of the velocity and magnetic field, we need the smallness condition on the $L^q$-norm of the gradient of the viscosity (see Lemma \ref{lem33}). In the special case that $\mu$ is a positive constant, it is clear that \eqref{C} holds true. Hence, Theorem \ref{thm1.1} implies that for any given (large) initial data $(\rho_0,\mathbf{u}_0,\mathbf{b}_0)$ satisfying \eqref{A}, there exists a unique global strong solution to the problem \eqref{1.1}--\eqref{1.3} with constant viscosity $\mu$ (see \cite{Z2021}).
\end{remark}


As a direct consequence of Theorem \ref{thm1.1}, we have the following global well-posedness and exponential decay of 2D nonhomogeneous Navier-Stokes equations with density-dependent viscosity.
\begin{theorem}\label{thm1.2}
Let the initial data $(\rho_0\geq0,\mathbf{u}_0)$ satisfy
\begin{align}\label{A2}
\rho_0\in W^{1,q},\
\mathbf{u}_0\in H_{0,\sigma}^{1}(\Omega),\ q\in(2,\infty).
\end{align}
Then there exists a small positive constant $\varepsilon_0$ depending only on $\Omega,\underline{\mu},\bar{\mu}:=\sup\limits_{[0,\|\rho_0\|_{L^\infty}]}
\mu(\rho), \nu, q, \|\rho_0\|_{L^\infty}$, and $\|\nabla\mathbf{u}_0\|_{L^2}^2$ such that if
\begin{align}\label{C2}
\|\nabla\mu(\rho_0)\|_{L^q}\leq\varepsilon_0,
\end{align}
the nonhomogeneous Navier-Stokes equations \eqref{1.1}--\eqref{1.3} with $\mathbf{b}=\mathbf{0}$ have a unique global strong solution $(\rho\geq0,\mathbf{u})$, which satisfies \eqref{1.5} and \eqref{lv1.2} with $\mathbf{b}=\mathbf{0}$ and $\sigma=\frac{\underline{\mu}}{d^2\|\rho_0\|_{L^\infty}}$.
\end{theorem}

\begin{remark}
Compared with \cite{HW2014}, on one hand, there is no need to impose the compatibility condition on the initial data despite the presence of vacuum. On the other hand, the exponential decay rate of the solution is a new result.
\end{remark}

\begin{remark}
We remark that the smallness condition \eqref{C2} allows
any given large initial data $(\rho_0,\mathbf{u}_0)$ satisfying \eqref{A2}
provided that the viscosity $\mu$ is a positive constant,
which is in sharp contrast to \cite{L2020} where the smallness assumption on $\|\rho_0\|_{L^\infty}$ is needed.
\end{remark}

\begin{remark}
It is not hard to prove that the strong-weak uniqueness theorem \cite[Theorem 2.7]{L1996} still holds for the initial data $(\rho_0,\mathbf{u}_0)$ satisfying \eqref{A2} after modifying its proof slightly. Therefore, our Theorem \ref{thm1.2} can be regarded as the uniqueness and regularity theory of Lions's weak solutions \cite{L1996} in 2D case
with $\nabla\mu(\rho_0)$ suitably small in the $L^q$-norm.
\end{remark}

We now make some comments on the key ingredients of the proof of Theorem \ref{thm1.1}. The local existence and uniqueness of strong solutions to the problem \eqref{1.1}--\eqref{1.3} follows from \cite{S2018} (see Lemma \ref{lem20}). Thus our efforts are devoted to establishing global a priori estimates on solutions in suitable higher-order norms. We will adapt some basic idea used in Huang and Wang \cite{HW2014}, where they investigated the global existence of strong solutions to the 2D nonhomogeneous Navier-Stokes equations with density-dependent viscosity and vacuum. However, compared with \cite{HW2014}, the proof of Theorem \ref{thm1.1} is much more involved due to the strong coupling between the velocity and the magnetic field and the absence of the compatibility condition \eqref{LHY}. Consequently, some new ideas are needed to overcome these difficulties.

As mentioned by \cite{HW2015,Z2015}, the key
ingredient here is to get the time-independent bounds on the $L^1(0,T;L^\infty)$-norm of $\nabla\mathbf{u}$ and then the $L^\infty(0,T;L^q)$-norm of $\nabla\mu(\rho)$.
First, applying the upper bounds on the density (see \eqref{3.3}) and the Poincar{\'e} inequality, we derive that $\|\sqrt{\rho}\mathbf{u}\|^2_{L^2}
+\|\mathbf{b}\|^2_{L^2}$ decays with the rate of $e^{-\sigma t}$ for some $\sigma>0$ depending only on $\underline{\mu},\nu,\|\rho_0\|_{L^\infty}$, and the diameter of the $\Omega$ (see \eqref{z3.3}). Next, we need to obtain time-weighted estimates of $\|\nabla\mathbf{u}\|_{L^2}^{2}+\|\nabla\mathbf{b}\|_{L^2}^{2}$.
To this end, we assume $\|\nabla\mu(\rho)\|_{L^q}\leq1$ on
$[0,T]$. Motivated by \cite{HW2014}, we make use of Desjardins' interpolation inequality (see Lemma \ref{lem24}) to control $\|\sqrt{\rho}\mathbf{u}\|_{L^4}$, and we find the key point is to control the term $\int\mathbf{b}\cdot\nabla\mathbf{b}\cdot\mathbf{u}_tdx$ (see \eqref{3.14}).
Multiplying \eqref{3.2}$_3$ by $\Delta\mathbf{b}$,
the term $\int\mathbf{b}\cdot\nabla\mathbf{b}\cdot\mathbf{u}_tdx$ can be controlled after using delicate Gagliardo-Nirenberg inequality (see \eqref{lv3.90}).
Next, using the structure of the 2D magnetic equation, we multiply \eqref{3.2}$_4$ by $|\mathbf{b}|^2\mathbf{b}$ and thus obtain some useful a priori estimates on $\||\mathbf{b}||\nabla\mathbf{b}|\|_{L^2}$, which is crucial in deriving the time-independent estimates on both the $L^\infty(0,T;L^2)$-norm of $t\sqrt{\rho}\mathbf{u}_t$ and the $L^2(0,T;L^2)$-norm of $t\nabla\mathbf{u}$ (see \eqref{4.16}). In fact, all these decay-in-time rates play an important role in obtaining the desired uniform bound (with respect to time) on the $L^1(0,T;L^\infty)$-norm of $\nabla\mathbf{u}$ (see \eqref{4.29}), which in particular implies $L^\infty(0,T;L^q)$-norm of the gradient of the viscosity $\mu(\rho)$ provided $\|\nabla\mu(\rho_0)\|_{L^q}\leq\varepsilon_0$ as stated in Theorem \ref{thm1.1} (see \eqref{3.21}). Finally, the higher order estimates on solutions are obtained (see Lemma \ref{lem35}) by considering time weighted type due to the lacking of the compatibility conditions.

The rest of this paper is organized as follows. In Section \ref{sec2}, we collect some elementary facts and inequalities that will be used later. Section \ref{sec3} is devoted to the
a priori estimates. Finally, we give the proof of Theorem \ref{thm1.1} in Section \ref{sec4}.

\section{Preliminaries}\label{sec2}

In this section, we will recall some known facts and elementary inequalities that will be used frequently later.

We begin with the local existence and uniqueness of strong solutions
whose proof can be found in \cite{S2018}.
\begin{lemma}\label{lem20}
Assume that $(\rho_0, \mathbf{u}_0,\mathbf{b}_0)$ satisfies \eqref{A}. Then there exist a small time $T>0$ and a unique strong solution $(\rho, \mathbf{u},\mathbf{b})$ to the problem \eqref{1.1}--\eqref{1.3} in $\Omega\times(0,T)$.
\end{lemma}

Next, the following Gagliardo-Nirenberg inequality (see \cite[Theorem 10.1, p. 27]{F2008}) will be useful in the next section.
\begin{lemma}[Gagliardo-Nirenberg]\label{lem22}
Let $\Omega\subset\mathbb{R}^2$ be a bounded smooth domain.
Assume that $1\leq q,r\leq\infty$, and $j,m$ are arbitrary integers satisfying $0\leq j<m$. If $v\in W^{m,r}(\Omega)\cap L^q(\Omega)$, then we have
\begin{equation*}
\|D^jv\|_{L^p}\leq C\|v\|_{L^q}^{1-a}\|v\|_{W^{m,r}}^{a},
\end{equation*}
where
\begin{equation*}
-j+\frac{2}{p}=(1-a)\frac{2}{q}+a\Big(-m+\frac{2}{r}\Big),
\end{equation*}
and
\begin{equation*}
\begin{split}
a\in
\begin{cases}
[\frac{j}{m},1),\ \ \text{if}\ m-j-\frac{2}{r}\ \text{is a nonnegative integer},\\
[\frac{j}{m},1],\ \ \text{otherwise}.
\end{cases}
\end{split}
\end{equation*}
The constant $C$ depends only on $m,j,q,r,a$, and $\Omega$. In particular, we have
\begin{align}\label{2.0}
\|v\|_{L^4}^4\leq C\|v\|_{L^2}^2\|v\|_{H^{1}}^2,
\end{align}
which will be used frequently in the next section.
\end{lemma}


Next, we give some regularity results for the following Stokes system with variable viscosity coefficient
\begin{align}\label{2.2}
\begin{cases}
-\divv(2\mu(\rho)\mathfrak{D}(\mathbf{u}))+\nabla P=\mathbf{F},\ \ x\in\Omega,\\
 \divv\mathbf{u}=0,\ \ x\in\Omega,\\
\mathbf{u}=\mathbf{0},\ \ x\in\partial\Omega, \\
\int Pdx=0.
\end{cases}
\end{align}
\begin{lemma}\label{lem23}
Assume that $\rho\in W^{1,q}(\Omega)$ with $2<q<\infty$, $0\leq\rho\leq\bar{\rho}$, $\mu\in C^1[0,\infty)$, and $\underline{\mu}\leq\mu(\rho)\leq\bar{\mu}$. Let $(\mathbf{u},P)\in H_{0}^1\times L^2$ be the unique weak solution to the problem \eqref{2.2}, then there exists a positive constant $C$ depending only on $\Omega,\bar{\rho},\underline{\mu},\bar{\mu}$ such that the following regularity results hold true:
\begin{itemize}
\item [$\bullet$] If $\mathbf{F}\in L^{2}(\Omega)$, then $(\mathbf{u}, P)\in H^2\times H^1$ and
\begin{equation*}
\begin{split}
& \|\mathbf{u}\|_{H^{2}}\leq C
\|\mathbf{F}\|_{L^2}\left(1+\|\nabla\mu(\rho)\|_{L^q}\right)^{\frac{q}{q-2}}, \\
& \|P\|_{H^{1}}\leq C
\|\mathbf{F}\|_{L^2}\left(1+\|\nabla\mu(\rho)\|_{L^q}\right)^{\frac{2q-2}{q-2}}.
\end{split}
\end{equation*}
\item [$\bullet$] If $\mathbf{F}\in L^{r}$ for some $r\in(2,q)$, then $(\mathbf{u},P)\in W^{2,r}\times W^{1,r}$ and
\begin{equation*}
\begin{split}
& \|\mathbf{u}\|_{W^{2,r}}
\leq C
\|\mathbf{F}\|_{L^r}\left(1+\|\nabla\mu(\rho)\|_{L^q}\right)^{\frac{qr}{2(q-r)}},
\\
& \|P\|_{W^{1,r}}
\leq C
\|\mathbf{F}\|_{L^r}
\left(1+\|\nabla\mu(\rho)\|_{L^q}\right)^{1+\frac{qr}{2(q-r)}}.
\end{split}
\end{equation*}
\end{itemize}
\end{lemma}
{\it Proof.}
See  \cite{CK2004,HW2014}.  \hfill $\Box$

Finally, by zero extension of $\mathbf{u}$ outside $\Omega$,
we can derive the following lemma  due to Desjardins
(see \cite[Lemma 1]{D1997} or \cite{D19972}), which plays a key role in the proof of Lemma \ref{lem33} in the next section.
\begin{lemma}\label{lem24}
Let $\Omega\subset\mathbb{R}^2$ be a bounded smooth domain.
Suppose that $0\leq\rho\leq\bar{\rho}$ and $\mathbf{u}\in H^1_0(\Omega)$, then we have
\begin{equation}\label{2.1}
\|\sqrt{\rho}\mathbf{u}\|_{L^4}^2\leq C(\bar{\rho},\Omega)(1+\|\sqrt{\rho}\mathbf{u}\|_{L^2})\|\nabla\mathbf{u}\|_{L^2}
\sqrt{\log(2+\|\nabla\mathbf{u}\|_{L^2}^2)}.
\end{equation}
\end{lemma}

\section{A priori estimates}\label{sec3}

In this section, we will establish some necessary a priori bounds for strong solutions $(\rho,\mathbf{u},\mathbf{b})$ to the problem \eqref{1.1}--\eqref{1.3} to extend the local strong solution. Thus, let $T>0$ be a fixed time and $(\rho,\mathbf{u},\mathbf{b})$  be the strong solution to \eqref{1.1}--\eqref{1.3} on $\Omega\times(0,T]$ with initial data $(\rho_0,\mathbf{u}_0,\mathbf{b}_0)$ satisfying \eqref{A}.
Before proceeding, we rewrite another equivalent form of the system \eqref{1.1} as the following
\begin{align}\label{3.2}
\begin{cases}
\rho_{t}+\mathbf{u}\cdot\nabla\rho=0,\\
\rho\mathbf{u}_{t}+\rho\mathbf{u}\cdot\nabla\mathbf{u}
-\divv(2\mu(\rho)\mathfrak{D}(\mathbf{u}))+\nabla P=\mathbf{b}\cdot\nabla\mathbf{b},\\
\mathbf{b}_t-\nu\Delta\mathbf{b}+\mathbf{u}\cdot\nabla\mathbf{b}
-\mathbf{b}\cdot\nabla\mathbf{b}=\mathbf{0},\\
\divv\mathbf{u}=\divv\mathbf{b}=0.
\end{cases}
\end{align}
In what follows, we denote by
\begin{equation*}
\int\cdot dx=\int_{\Omega}\cdot dx.
\end{equation*}
We sometimes use $C(f)$ to emphasize the dependence on $f$.

First, since \eqref{3.2}$_1$ is a transport equation, we have directly the following result.
\begin{lemma}\label{lem31}
For $(x,t)\in\Omega\times[0,T]$, it holds that
\begin{align}\label{3.3}
0\leq \rho(x,t)\leq \|\rho_0\|_{L^\infty}.
\end{align}
\end{lemma}
\begin{remark}\label{re1}
Since $\mu(\rho)$ is a continuously differentiable function, we deduce from \eqref{3.3} and \eqref{1.4} that
\begin{equation}\label{z3.2}
0<\underline{\mu}\leq\mu(\rho)\leq\bar{\mu}
:=\sup_{[0,\|\rho_0\|_{L^\infty}]}\mu(\rho)<\infty,
\end{equation}
and
\begin{equation*}
\|\mu'(\rho)\|_{L^\infty(0,T;L^\infty)}<\infty.
\end{equation*}
\end{remark}

Next, the following lemma gives the basic energy estimates.
\begin{lemma}\label{lem32}
It holds that
\begin{equation}\label{3.4}
\sup_{0\leq t\leq T}\left(\|\sqrt{\rho}\mathbf{u}\|_{L^2}^2
+\|\mathbf{b}\|_{L^2}^2\right)
+\int_{0}^{T}\big(\underline{\mu}\|\nabla\mathbf{u}\|_{L^2}^2
+\nu\|\nabla\mathbf{b}\|_{L^2}^2\big)dt
\leq \|\sqrt{\rho_0}\mathbf{u}_0\|_{L^2}^2
+\|\mathbf{b}_0\|_{L^2}^2,
\end{equation}
and
\begin{align}\label{z3.3}
\sup_{0\leq t\leq T}\big[e^{\sigma t}
\left(\|\sqrt{\rho}\mathbf{u}\|_{L^2}^2
+\|\mathbf{b}\|_{L^2}^2\right)\big]
+\int_{0}^{T}e^{\sigma t}\big(\underline{\mu}\|\nabla\mathbf{u}\|_{L^2}^2
+\nu\|\nabla\mathbf{b}\|_{L^2}^2\big)dt
\leq \|\sqrt{\rho_0}\mathbf{u}_0\|_{L^2}^2+\|\mathbf{b}_0\|_{L^2}^2,
\end{align}
where $\sigma:=\min\left\{\frac{\underline{\mu}}{d^2\|\rho_0\|_{L^\infty}},
\frac{\nu}{d^2}\right\}$ with $d$ being the diameter of $\Omega$.
\end{lemma}
{\it Proof.}
1. Multiplying \eqref{3.2}$_2$ by $\mathbf{u}$, \eqref{3.2}$_3$ by $\mathbf{b}$, and integration (by parts) over $\Omega$, we derive that
\begin{equation}\label{3.5}
\frac12\frac{d}{dt}\left(\|\sqrt{\rho}\mathbf{u}\|_{L^2}^2
+\|\mathbf{b}\|_{L^2}^2\right)
+2\int\mu(\rho)\mathfrak{D}(\mathbf{u})\cdot\nabla\mathbf{u}dx
+\nu\|\nabla\mathbf{b}\|_{L^2}^2=0.
\end{equation}
Noting that
\begin{align*}
2\int\mu(\rho)\mathfrak{D}(\mathbf{u})\cdot\nabla\mathbf{u}dx
& =\int\mu(\rho)(\partial_iu^j+\partial_ju^i)\partial_iu^jdx \nonumber \\
& = \frac{1}{2}\int\mu(\rho)(\partial_iu^j+\partial_ju^i)
(\partial_iu^j+\partial_ju^i)dx
\nonumber \\
& = 2\int\mu(\rho)|\mathfrak{D}(\mathbf{u})|^2dx,
\end{align*}
and
\begin{align}\label{z3.4}
2\int|\mathfrak{D}(\mathbf{u})|^2dx
& =\frac{1}{2}\int(\partial_iu^j+\partial_ju^i)(\partial_iu^j+\partial_ju^i)dx \notag \\
& =\int|\nabla\mathbf{u}|^2dx+\int\partial_iu^j\partial_ju^idx
=\int|\nabla\mathbf{u}|^2dx,
\end{align}
we thus obtain from \eqref{3.5} and \eqref{1.4} that
\begin{equation}\label{3.6}
\frac{d}{dt}\left(\|\sqrt{\rho}\mathbf{u}\|_{L^2}^2
+\|\mathbf{b}\|_{L^2}^2\right)
+2\big(\underline{\mu}\|\nabla\mathbf{u}\|_{L^2}^2
+\nu\|\nabla\mathbf{b}\|_{L^2}^2\big)\leq0.
\end{equation}
Integrating the above inequality over $(0,T)$ gives the desired \eqref{3.4}.

2. It follows from Poincar{\'e}'s inequality (see \cite[(A.3), p. 266]{S2008}) and \eqref{3.3} that
\begin{align}\label{zz3.1}
\|\sqrt{\rho}\mathbf{u}\|_{L^2}^2
\leq\|\rho\|_{L^\infty}\|\mathbf{u}\|_{L^2}^2
\leq\|\rho_0\|_{L^\infty}d^2\|\nabla\mathbf{u}\|_{L^2}^2,
\end{align}
where $d$ is the diameter of $\Omega$. Hence, we get
\begin{align}\label{zx1}
\frac{1}{d^2\|\rho_0\|_{L^\infty}}\|\sqrt{\rho}\mathbf{u}\|_{L^2}^2
\leq\|\nabla\mathbf{u}\|_{L^2}^2,\ \frac{1}{d^2}\|\mathbf{b}\|_{L^2}^2
\leq \|\nabla\mathbf{b}\|_{L^2}^2.
\end{align}
Consequently, letting $\sigma:=\min\left\{\frac{\underline{\mu}}{d^2\|\rho_0\|_{L^\infty}},
\frac{\nu}{d^2}\right\}$, then we derive from \eqref{3.6} and \eqref{zx1} that
\begin{equation}\label{3.7}
\frac{d}{dt}\big[e^{\sigma t}\left(\|\sqrt{\rho}\mathbf{u}\|_{L^2}^2
+\|\mathbf{b}\|_{L^2}^2\right)\big]
+e^{\sigma t}\left(\underline{\mu}\|\nabla\mathbf{u}\|_{L^2}^2
+\nu\|\nabla\mathbf{b}\|_{L^2}^2\right)
\leq0.
\end{equation}
Thus, integrating \eqref{3.7} with respect to $t$ gives \eqref{z3.3}.
\hfill $\Box$

\begin{lemma}\label{lem33}
Let $q$ be as in Theorem \ref{thm1.1} and $\bar{\mu}$ be as in \eqref{z3.2}, assume that
\begin{align}\label{3.1}
\sup_{0\le t\le T}\|\nabla\mu(\rho)\|_{L^q}\leq 1,
\end{align}
then there exists a positive constant $C$ depending only on $\Omega$, $\underline{\mu}$, $\bar{\mu}$, $\nu$, $q$, $\|\rho_0\|_{L^\infty}$, $\|\nabla\mathbf{u}_0\|_{L^2}^2$, and $\|\nabla\mathbf{b}_0\|_{L^2}^2$ such that
\begin{align}\label{3.9}
\sup_{0\leq t\leq T}\left(\|\nabla\mathbf{u}\|_{L^2}^{2}
+\|\nabla\mathbf{b}\|_{L^2}^{2}\right)
+\int_{0}^{T}\left(\|\sqrt{\rho}\mathbf{u}_t\|_{L^2}^{2}
+\|\Delta\mathbf{b}\|_{L^2}^{2}
+\||\mathbf{b}||\nabla\mathbf{b}|\|_{L^2}^2\right)dt \leq C.
\end{align}
Moreover, for $\sigma$ as in Lemma \ref{lem32}, one has
\begin{align}\label{z3.10}
\sup_{0\leq t\leq T}\big[e^{\sigma t}\left(\|\nabla\mathbf{u}\|_{L^2}^{2}
+\|\nabla\mathbf{b}\|_{L^2}^{2}\right)\big]
+\int_{0}^{T}e^{\sigma t}\left(\|\sqrt{\rho}\mathbf{u}_t\|_{L^2}^{2}
+\|\Delta\mathbf{b}\|_{L^2}^{2}
+\||\mathbf{b}||\nabla\mathbf{b}|\|_{L^2}^2\right)dt\leq C.
\end{align}
\end{lemma}
{\it Proof.}
1. Since $\mu(\rho)$ is a continuously differentiable
function, we obtain from \eqref{3.2}$_1$ that
\begin{equation}\label{3.10}
[\mu(\rho)]_t+\mathbf{u}\cdot\nabla\mu(\rho)=0.
\end{equation}
Multiplying \eqref{3.2}$_2$ by $\mathbf{u}_{t}$ and integrating the resulting equation over $\Omega$ imply that
\begin{equation}\label{3.11}
2\int\mu(\rho)\mathfrak{D}(\mathbf{u})\cdot\nabla\mathbf{u}_tdx
+\int\rho|\mathbf{u}_{t}|^2dx
= -\int\rho\mathbf{u}\cdot\nabla\mathbf{u}\cdot\mathbf{u}_{t}dx
+\int\mathbf{b}\cdot\nabla\mathbf{b}\cdot\mathbf{u}_tdx.
\end{equation}
Similarly to \eqref{3.7}, we get
\begin{equation*}
\begin{split}
2\int\mu(\rho)\mathfrak{D}(\mathbf{u})\cdot\nabla\mathbf{u}_tdx
 =2\int\mu(\rho)\mathfrak{D}(\mathbf{u})\cdot\mathfrak{D}(\mathbf{u})_tdx
 =\frac{d}{dt}\int\mu(\rho)|\mathfrak{D}(\mathbf{u})|^2dx
-\int[\mu(\rho)]_t|\mathfrak{D}(\mathbf{u})|^2dx,
\end{split}
\end{equation*}
which combined with \eqref{3.11} and \eqref{3.10} leads to
\begin{equation}\label{3.12}
\frac{d}{dt}\int\mu(\rho)|\mathfrak{D}(\mathbf{u})|^2dx
+\int\rho|\mathbf{u}_{t}|^2dx
= -\int\rho\mathbf{u}\cdot\nabla\mathbf{u}\cdot\mathbf{u}_{t}dx
-\int\mathbf{u}\cdot\nabla\mu(\rho)|\mathfrak{D}(\mathbf{u})|^2dx
+\int\mathbf{b}\cdot\nabla\mathbf{b}\cdot\mathbf{u}_tdx.
\end{equation}
By H{\"o}lder's and Gagliardo-Nirenberg inequalities, we have
\begin{align}\label{3.13}
\left|-\int\rho\mathbf{u}\cdot\nabla\mathbf{u}\cdot\mathbf{u}_{t}dx\right|
& \leq \frac{1}{4}\|\sqrt{\rho}\mathbf{u}_t\|_{L^2}^2
+2\|\sqrt{\rho}\mathbf{u}\|_{L^4}^2\|\nabla\mathbf{u}\|_{L^4}^2 \nonumber \\
& \leq \frac{1}{4}\|\sqrt{\rho}\mathbf{u}_t\|_{L^2}^2
+C(\Omega)\|\sqrt{\rho}\mathbf{u}\|_{L^4}^2
\|\nabla\mathbf{u}\|_{L^2}\|\mathbf{u}\|_{H^2}.
\end{align}
By Sobolev's inequality, \eqref{3.1}, and Gagliardo-Nirenberg inequality, we arrive at
\begin{align}
\left|-\int\mathbf{u}\cdot\nabla\mu(\rho)|\mathfrak{D}(\mathbf{u})|^2dx\right|
& \leq C\int|\mathbf{u}||\nabla\mu(\rho)||\nabla\mathbf{u}|^2dx \nonumber \\
& \leq C\|\nabla\mu(\rho)\|_{L^q}\|\mathbf{u}\|_{L^{\frac{2q}{q-2}}}
\|\nabla\mathbf{u}\|_{L^4}^2
\nonumber \\
& \leq C(\Omega)
\|\nabla\mathbf{u}\|_{L^2}\|\nabla\mathbf{u}\|_{L^2}
\|\nabla\mathbf{u}\|_{H^1}
\nonumber \\
& \leq C\|\nabla\mathbf{u}\|_{L^2}^2\|\mathbf{u}\|_{H^2}.
\end{align}
Integration by parts together with $\divv\mathbf{b}=0$ and $\mathbf{b}|_{\partial\Omega}=\mathbf{0}$, we infer from Sobolev's inequality that
\begin{align}\label{3.14}
\int\mathbf{b}\cdot\nabla\mathbf{b}\cdot\mathbf{u}_tdx
& = -\frac{d}{dt}\int \mathbf{b}\cdot\nabla\mathbf{u}\cdot\mathbf{b}dx
+\int\mathbf{b}_t\cdot\nabla\mathbf{u}\cdot\mathbf{b}dx
+\int\mathbf{b}\cdot\nabla\mathbf{u}\cdot\mathbf{b}_tdx \notag \\
& = -\frac{d}{dt}\int \mathbf{b}\cdot\nabla\mathbf{u}\cdot\mathbf{b}dx
+\int(\nu\Delta\mathbf{b}-\mathbf{u}\cdot\nabla\mathbf{b}+\mathbf{b}\cdot\nabla \mathbf{u})\cdot\nabla\mathbf{u}\cdot\mathbf{b}dx \notag \\
& \quad
+\int\mathbf{b}\cdot\nabla\mathbf{u}\cdot
(\nu\Delta\mathbf{b}-\mathbf{u}\cdot\nabla\mathbf{b}+\mathbf{b}\cdot\nabla \mathbf{u})dx \notag \\
& \leq -\frac{d}{dt}\int \mathbf{b}\cdot\nabla\mathbf{u}\cdot\mathbf{b}dx
+\frac{\nu}{4}\|\Delta\mathbf{b}\|_{L^2}^2+C\|\mathbf{b}\|_{L^6}^6
+C\|\nabla\mathbf{u}\|_{L^3}^{3} \notag \\
& \quad +C\|\mathbf{b}\|_{L^4}\|\mathbf{u}\|_{L^\infty}
\|\nabla\mathbf{b}\|_{L^2}\|\nabla\mathbf{u}\|_{L^4}
\notag \\
& \leq -\frac{d}{dt}\int \mathbf{b}\cdot\nabla\mathbf{u}\cdot\mathbf{b}dx
+\frac{\nu}{4}\|\Delta\mathbf{b}\|_{L^2}^2
\notag \\
& \quad +C\|\mathbf{b}\|_{L^2}^2\|\nabla\mathbf{b}\|_{L^2}^4
+C\|\nabla\mathbf{u}\|_{L^2}^{2}\|\mathbf{u}\|_{H^2} \notag \\
& \quad
+C\|\mathbf{b}\|_{L^4}\|\mathbf{u}\|_{L^4}^{\frac12}
\|\nabla\mathbf{u}\|_{L^4}^{\frac12}
\|\nabla\mathbf{b}\|_{L^2}\|\nabla\mathbf{u}\|_{L^2}^{\frac12}
\|\nabla\mathbf{u}\|_{H^1}^{\frac12} \notag \\
& \leq -\frac{d}{dt}\int \mathbf{b}\cdot\nabla\mathbf{u}\cdot\mathbf{b}dx
+\frac{\nu}{4}\|\Delta\mathbf{b}\|_{L^2}^2+C\|\mathbf{b}\|_{L^2}^2\|\nabla\mathbf{b}\|_{L^2}^4
\notag \\
& \quad +C\|\nabla\mathbf{u}\|_{L^2}^{2}\|\mathbf{u}\|_{H^2} +C\|\mathbf{b}\|_{L^4}
\|\nabla\mathbf{b}\|_{L^2}\|\nabla\mathbf{u}\|_{L^2}
\|\mathbf{u}\|_{H^2},
\end{align}
where we have used the following Gagliardo-Nirenberg inequality
\begin{align*}
\|\mathbf{u}\|_{L^\infty}\leq C\|\mathbf{u}\|_{L^4}^{\frac12}
\|\nabla\mathbf{u}\|_{L^4}^{\frac12},\
\|\nabla\mathbf{u}\|_{L^4}\leq C\|\nabla\mathbf{u}\|_{L^2}^{\frac12}
\|\nabla\mathbf{u}\|_{H^1}^{\frac12}.
\end{align*}
Substituting \eqref{3.13}--\eqref{3.14} into \eqref{3.12}, we derive
\begin{align}\label{3.15}
& \frac{d}{dt}\int\big(\mu(\rho)|\mathfrak{D}(\mathbf{u})|^2
+\mathbf{b}\cdot\nabla\mathbf{u}\cdot\mathbf{b}\big)dx
+\frac12\|\sqrt{\rho}\mathbf{u}_{t}\|_{L^2}^2 \nonumber \\
& \leq C\left(\|\sqrt{\rho}\mathbf{u}\|_{L^4}^2+\|\nabla\mathbf{u}\|_{L^2}\right)
\|\nabla\mathbf{u}\|_{L^2}\|\mathbf{u}\|_{H^2}
+\frac{\nu}{4}\|\Delta\mathbf{b}\|_{L^2}^2+C\|\nabla\mathbf{b}\|_{L^2}^4 \notag \\
& \quad +C\|\mathbf{b}\|_{L^4}
\|\nabla\mathbf{b}\|_{L^2}\|\nabla\mathbf{u}\|_{L^2}
\|\mathbf{u}\|_{H^2}.
\end{align}

2. Multiplying \eqref{3.2}$_3$ by $\Delta\mathbf{b}$ and integrating the resulting equality over $\Omega$,
it follows from H\"older's and Gagliardo-Nirenberg inequalities that
\begin{align}\label{mm1}
& \frac{d}{dt}\int|\nabla\mathbf{b}|^2dx
+\nu\int|\Delta\mathbf{b}|^2dx \notag \\
& \leq C\int |\nabla\mathbf{u}||\nabla\mathbf{b}|^2dx
+C\int|\nabla\mathbf{u}||\mathbf{b}||\Delta\mathbf{b}|dx \notag \\
& \leq C\|\nabla\mathbf{u}\|_{L^3}\|\nabla\mathbf{b}\|_{L^2}^{\frac43}
\|\Delta\mathbf{b}\|_{L^2}^{\frac23}+C\|\nabla\mathbf{u}\|_{L^3}\|\mathbf{b}\|_{L^6}
\|\Delta\mathbf{b}\|_{L^2}\notag \\
& \leq C\|\nabla\mathbf{u}\|_{L^2}^2\|\nabla^2\mathbf{u}\|_{L^2}+C (1+\|\mathbf{b}\|_{L^2}^2)\|\nabla\mathbf{b}\|_{L^2}^4
+\frac{\nu}{4}\|\Delta\mathbf{b}\|_{L^2}^{2},
\end{align}
which together with \eqref{3.15} and \eqref{3.4} gives rise to
\begin{align}\label{lv3.90}
& B'(t)+\frac12\|\sqrt{\rho}\mathbf{u}_t\|_{L^2}^{2}
+\frac{\nu}{2}\|\Delta\mathbf{b}\|_{L^2}^{2} \notag \\
& \leq C\left(\|\sqrt{\rho}\mathbf{u}\|_{L^4}^2+\|\nabla\mathbf{u}\|_{L^2}\right)
\|\nabla\mathbf{u}\|_{L^2}\|\mathbf{u}\|_{H^2}
+C\|\nabla\mathbf{b}\|_{L^2}^4+C\|\mathbf{b}\|_{L^4}
\|\nabla\mathbf{b}\|_{L^2}\|\nabla\mathbf{u}\|_{L^2}
\|\mathbf{u}\|_{H^2},
\end{align}
where
\begin{equation}\label{lv3.11}
B(t):= \int\big(\mu(\rho)|\mathfrak{D}(\mathbf{u})|^2
+|\nabla\mathbf{b}|^2+\mathbf{b}\cdot\nabla\mathbf{u}\cdot\mathbf{b}\big)dx
\end{equation}
satisfies
\begin{equation}\label{3.20}
\frac{\underline{\mu}}{4}\|\nabla\mathbf{u}\|_{L^2}^{2}+\|\nabla\mathbf{b}\|_{L^2}^{2}
-C_1\|\mathbf{b}\|_{L^4}^{4}
\leq B(t)
\leq C\|\nabla\mathbf{u}\|_{L^2}^{2}+C\|\nabla\mathbf{b}\|_{L^2}^{2}
\end{equation}
owing to Gagliardo-Nirenberg inequality, \eqref{3.4}, and the following estimate
\begin{equation}\label{llv3.11}
\int |\mathbf{b}\cdot\nabla\mathbf{u}\cdot\mathbf{b}|dx\le  \frac{\underline{\mu}}{4}\|\nabla\mathbf{u}\|_{L^2}^{2}+C_1\|\mathbf{b}\|_{L^4}^{4}.
\end{equation}

3. Multiplying\eqref{3.2}$_3$ by $|\mathbf{b}|^2\mathbf{b}$ and integrating the resulting equality by parts over $\Omega$, we obtain from Gagliardo-Nirenberg inequality that
\begin{align}\label{lv3.7}
\frac{1}{4}\frac{d}{dt}\|\mathbf{b}\|^4_{L^4}+\||\nabla\mathbf{b}| |\mathbf{b}|\|_{L^2}^2+\frac{1}{2}\|\nabla|\mathbf{b}|^2\|_{L^2}^2
& \leq C\|\nabla\mathbf{u}\|_{L^2} \||\mathbf{b}|^2\|_{L^4}^2 \notag \\
&\leq C\|\nabla\mathbf{u}\|_{L^2} \||\mathbf{b}|^2\|_{L^2}
\|\nabla|\mathbf{b}|^2\|_{L^2} \notag \\
&\leq \frac{1}{4}\|\nabla|\mathbf{b}|^2\|_{L^2}^2
+C\|\nabla\mathbf{u}\|_{L^2}^2\|\mathbf{b}\|_{L^4}^4,
\end{align}
which together with Gronwall's inequality and \eqref{3.4} implies
\begin{align}\label{lz}
& \sup_{0\leq t\leq T}\|\mathbf{b}\|_{L^4}^4
+\int_{0}^{T}\||\mathbf{b}||\nabla\mathbf{b}|\|_{L^2}^2dt \leq C.
\end{align}
This along with \eqref{lv3.90} yields
\begin{align}\label{lv3.9}
& B'(t)+\frac12\|\sqrt{\rho}\mathbf{u}_t\|_{L^2}^{2}
+\frac{\nu}{2}\|\Delta\mathbf{b}\|_{L^2}^{2} \notag \\
& \leq C\left(\|\sqrt{\rho}\mathbf{u}\|_{L^4}^2+\|\nabla\mathbf{u}\|_{L^2}\right)
\|\nabla\mathbf{u}\|_{L^2}\|\mathbf{u}\|_{H^2}
+C\|\nabla\mathbf{b}\|_{L^2}^4+C\|\nabla\mathbf{b}\|_{L^2}\|\nabla\mathbf{u}\|_{L^2}
\|\mathbf{u}\|_{H^2}.
\end{align}

4. Recall that $(\mathbf{u}, P)$ satisfies the following Stokes system with variable viscosity
\begin{equation*}
\begin{cases}
 -\divv(2\mu(\rho)\mathfrak{D}(\mathbf{u})) + \nabla P = -\rho\mathbf{u}_t-\rho\mathbf{u}\cdot\nabla\mathbf{u}
 +\mathbf{b}\cdot\nabla\mathbf{b},\,\,\,\,&x\in \Omega,\\
 \divv\mathbf{u}=0,   \,\,\,&x\in \Omega,\\
\mathbf{u}=\mathbf{0},\,\,\,\,&x\in \partial\Omega.
\end{cases}
\end{equation*}
Applying Lemma \ref{lem23} with $\mathbf{F}=-\rho\mathbf{u}_t-\rho\mathbf{u}\cdot\nabla\mathbf{u}
+\mathbf{b}\cdot\nabla\mathbf{b}$, we obtain from \eqref{3.3} and \eqref{3.1} that
\begin{equation*}
\begin{split}
\|\mathbf{u}\|_{H^2}+\|\nabla P\|_{L^2}
& \leq C\left(\|\rho\mathbf{u}_t\|_{L^2}
+\|\rho\mathbf{u}\cdot\nabla\mathbf{u}\|_{L^2}
+\|\mathbf{b}\cdot\nabla\mathbf{b}\|_{L^2}\right)
(1+\|\nabla\mu(\rho)\|_{L^q})^{\frac{q}{q-2}} \nonumber \\
& \leq C\|\sqrt{\rho}\mathbf{u}_t\|_{L^2}
+C\|\sqrt{\rho}\mathbf{u}\|_{L^4}\|\nabla\mathbf{u}\|_{L^4}
+C\||\mathbf{b}||\nabla\mathbf{b}|\|_{L^2}
 \nonumber \\
& \leq C\|\sqrt{\rho}\mathbf{u}_t\|_{L^2}
+C\|\sqrt{\rho}\mathbf{u}\|_{L^4}\|\nabla\mathbf{u}\|_{L^2}^{\frac{1}{2}}
\|\mathbf{u}\|_{H^2}^{\frac12}+C\||\mathbf{b}||\nabla\mathbf{b}|\|_{L^2}
 \nonumber \\
& \leq C\|\sqrt{\rho}\mathbf{u}_t\|_{L^2}
+C\|\sqrt{\rho}\mathbf{u}\|_{L^4}^2\|\nabla\mathbf{u}\|_{L^2}
+\frac{1}{2}\|\mathbf{u}\|_{H^2}+C\||\mathbf{b}||\nabla\mathbf{b}|\|_{L^2},
\end{split}
\end{equation*}
and thus
\begin{equation}\label{3.16}
\|\mathbf{u}\|_{H^2}+\|\nabla P\|_{L^2}
\leq C\|\sqrt{\rho}\mathbf{u}_t\|_{L^2}
+C\|\sqrt{\rho}\mathbf{u}\|_{L^4}^2\|\nabla\mathbf{u}\|_{L^2}
+C\||\mathbf{b}||\nabla\mathbf{b}|\|_{L^2}.
\end{equation}
Inserting \eqref{3.16} into \eqref{3.15} and applying Cauchy-Schwarz inequality, we deduce that
\begin{align}\label{p}
& B'(t)
+\frac12\|\sqrt{\rho}\mathbf{u}_{t}\|_{L^2}^2
+\frac{\nu}{2}\|\Delta\mathbf{b}\|_{L^2}^2 \notag \\
& \leq \frac14\|\sqrt{\rho}\mathbf{u}_t\|_{L^2}^2
+C\|\sqrt{\rho}\mathbf{u}\|_{L^4}^4\|\nabla\mathbf{u}\|_{L^2}^2
+C\|\nabla\mathbf{u}\|_{L^2}^4
+C\|\nabla\mathbf{b}\|_{L^2}^4
+\varepsilon\||\mathbf{b}||\nabla\mathbf{b}|\|_{L^2}^2.
\end{align}
Noting that
\begin{align*}
\|\nabla\mathbf{u}\|_{L^2}^2\|\mathbf{b}\|_{L^4}^4
\leq C\|\nabla\mathbf{u}\|_{L^2}^2\|\mathbf{b}\|_{L^2}^2\|\nabla\mathbf{b}\|_{L^2}^2
\leq C\|\nabla\mathbf{u}\|_{L^2}^4+C\|\nabla\mathbf{b}\|_{L^2}^4
\end{align*}
due to \eqref{3.4} and Cauchy-Schwarz inequality. Thus,
adding \eqref{lv3.7} multiplied by $4(C_1+1)$ to \eqref{p} and choosing $\varepsilon$ suitably small, we obtain after using \eqref{2.1} and \eqref{3.4} that
\begin{align}\label{3.17}
&\frac{d}{dt}\left(B(t)+ (C_1+1)\|\mathbf{b}\|_{L^4}^{4}\right)
+\|\sqrt{\rho}\mathbf{u}_t\|_{L^2}^{2}
+\nu\|\Delta\mathbf{b}\|_{L^2}^{2}
+\||\mathbf{b}||\nabla\mathbf{b}|\|_{L^2}^2 \notag \\
& \leq C\|\nabla\mathbf{b}\|_{L^2}^4+C\|\nabla\mathbf{u}\|_{L^2}^{4}
+C\|\sqrt{\rho}\mathbf{u}\|_{L^4}^4\|\nabla\mathbf{u}\|_{L^2}^{2}
\notag \\
& \leq C\|\nabla\mathbf{b}\|_{L^2}^2\|\nabla\mathbf{b}\|_{L^2}^2
+C\|\nabla\mathbf{u}\|_{L^2}^2\|\nabla\mathbf{u}\|_{L^2}^2 \notag \\
& \quad +C\|\nabla\mathbf{u}\|_{L^2}^2\|\nabla\mathbf{u}\|_{L^2}^2
\log(2+\|\nabla\mathbf{u}\|_{L^2}^2).
\end{align}
Set
\begin{equation*}
f(t):=2+B(t)+ (C_1+1)\|\mathbf{b}\|_{L^4}^{4},\
g(t):=\|\nabla\mathbf{u}\|_{L^2}^{2}+\|\nabla\mathbf{b}\|_{L^2}^{2},
\end{equation*}
then we deduce from \eqref{3.17} and \eqref{3.20} that
\begin{equation*}
f'(t)
\leq Cg(t)f(t)+Cg(t)f(t)\log f(t),
\end{equation*}
which yields
\begin{equation}\label{3.30}
(\log f(t))'\leq Cg(t)+Cg(t)\log(f(t)).
\end{equation}
We thus infer from \eqref{3.30}, Gronwall's inequality, \eqref{3.4}, and \eqref{3.20} that
\begin{equation}\label{3.31}
\sup_{0\leq t\leq T}\left(\|\nabla\mathbf{u}\|_{L^2}^{2}
+\|\nabla\mathbf{b}\|_{L^2}^{2}+\|\mathbf{b}\|_{L^4}^{4}\right)
 \leq C.
\end{equation}
Integrating \eqref{3.17} with respect to $t$ together with \eqref{3.31} and \eqref{3.4} leads to
\begin{equation}\label{3.32}
\int_{0}^{T}\left(\|\sqrt{\rho}\mathbf{u}_t\|_{L^2}^{2}
+\|\Delta\mathbf{b}\|_{L^2}^{2}
+\||\mathbf{b}||\nabla\mathbf{b}|\|_{L^2}^2\right)dt \leq C.
\end{equation}
This along with \eqref{3.31} gives the desired \eqref{3.9}.

5. Multiplying \eqref{3.17} by $e^{\sigma t}$ and applying \eqref{3.31}, we derive from \eqref{3.20} that
\begin{align}\label{z3.17}
&\frac{d}{dt}\left(e^{\sigma t}B(t)+ (C_1+1)e^{\sigma t}\|\mathbf{b}\|_{L^4}^{4}\right)
+e^{\sigma t}\left(\|\sqrt{\rho}\mathbf{u}_t\|_{L^2}^{2}
+\nu\|\Delta\mathbf{b}\|_{L^2}^{2}
+\||\mathbf{b}||\nabla\mathbf{b}|\|_{L^2}^2\right) \notag \\
& \leq Ce^{\sigma t}\big(\|\nabla\mathbf{b}\|_{L^2}^2
+\|\nabla\mathbf{u}\|_{L^2}^2\big)+\sigma e^{\sigma t}B(t)+ \sigma(C_1+1)e^{\sigma t}\|\mathbf{b}\|_{L^4}^{4} \notag \\
& \leq Ce^{\sigma t}\big(\|\nabla\mathbf{b}\|_{L^2}^2
+\|\nabla\mathbf{u}\|_{L^2}^2\big)+Ce^{\sigma t}\|\mathbf{b}\|_{L^2}^2
\|\nabla\mathbf{b}\|_{L^2}^2 \notag \\
& \leq Ce^{\sigma t}\big(\|\nabla\mathbf{b}\|_{L^2}^2
+\|\nabla\mathbf{u}\|_{L^2}^2\big).
\end{align}
Integrating \eqref{z3.17} over $(0,T)$ together with \eqref{3.20} leads to \eqref{z3.10}.
\hfill $\Box$

\begin{remark}\label{re2}
Under the condition \eqref{3.1}, it follows from \eqref{2.1}, \eqref{3.3}, \eqref{3.4}, and \eqref{3.9} that
\begin{align}\label{z3.9}
\sup_{0\leq t\leq T}\|\sqrt{\rho}\mathbf{u}\|_{L^4}^{2}\leq C.
\end{align}
\end{remark}

\begin{lemma}\label{lem43}
Let the condition \eqref{3.1} be satisfied, then there exists a positive constant $C$ depending only on $\Omega,\underline{\mu},\bar{\mu}, \nu, q, \|\rho_0\|_{L^\infty}, \|\nabla\mathbf{u}_0\|_{L^2}^2$, and $\|\nabla\mathbf{b}_0\|_{L^2}^2$ such that
for $i\in\{1,2\}$,
\begin{align}\label{4.16}
\sup_{0\leq t\leq T}\big[t^i\big(\|\sqrt{\rho}\mathbf{u}_t\|_{L^2}^2
+\|\mathbf{b}_t\|_{L^2}^2\big)\big]
+\int_{0}^Tt^i\big(\|\nabla\mathbf{u}_t\|_{L^2}^2
+\|\nabla\mathbf{b}_t\|_{L^2}^2\big)dt
\leq C.
\end{align}
Moreover, for $\sigma$ as that in Lemma \ref{lem32}, one has
\begin{align}\label{xz4.16}
\sup_{\zeta(T)\leq t\leq T}\big[e^{\sigma t}
\left(\|\sqrt{\rho}\mathbf{u}_t\|_{L^2}^2+\|\mathbf{b}_t\|_{L^2}^2\right)\big]
+\int_{\zeta(T)}^{T}e^{\sigma t}\left(\|\nabla\mathbf{u}_t\|_{L^2}^{2}
+\|\nabla\mathbf{b}_t\|_{L^2}^{2}\right)dt\leq C,
\end{align}
where $\zeta(T):=\min\{1,T\}$.
\end{lemma}
{\it Proof.}
1. Differentiating \eqref{3.2}$_2$ with respect to $t$, we arrive at
\begin{align}\label{3.22}
& \rho\mathbf{u}_{tt}+\rho\mathbf{u}\cdot\nabla\mathbf{u}_{t}
-\divv(2\mu(\rho)\mathfrak{D}(\mathbf{u}_t))\notag \\
& = -\nabla P_{t}
+\rho_t\left(\mathbf{u}_{t}+\mathbf{u}\cdot\nabla\mathbf{u}\right)
-\rho\mathbf{u}_{t}\cdot\nabla\mathbf{u}+\divv(2\mu_t\mathfrak{D}(\mathbf{u}))
+\mathbf{b}_{t}\cdot\nabla\mathbf{b}
+\mathbf{b}\cdot\nabla\mathbf{b}_{t}.
\end{align}
Multiplying \eqref{3.22} by $\mathbf{u}_t$ and integrating (by parts) over $\Omega$  and using \eqref{1.1}$_1$ yield
\begin{align}\label{3.23}
& \frac{1}{2}\frac{d}{dt}\int\rho|\mathbf{u}_{t}|^2dx
+2\int\mu(\rho)\mathfrak{D}(\mathbf{u}_t)\cdot\nabla\mathbf{u}_{t}dx
\nonumber \\ &
=\int\divv(\rho\mathbf{u})|\mathbf{u}_{t}|^2dx
+\int\divv(\rho\mathbf{u})\mathbf{u}\cdot\nabla\mathbf{u}\cdot\mathbf{u}_{t}dx
-\int\rho\mathbf{u}_{t}\cdot\nabla\mathbf{u}\cdot\mathbf{u}_{t}dx
-\int2\mu_t\mathfrak{D}(\mathbf{u})\cdot\nabla\mathbf{u}_tdx\nonumber \\
& \quad +\int\mathbf{b}_{t}\cdot\nabla\mathbf{b}\cdot\mathbf{u}_{t}dx
+\int\mathbf{b}\cdot\nabla\mathbf{b}_{t}\cdot\mathbf{u}_{t}dx
=:\sum_{i=1}^6J_i.
\end{align}
By virtue of H{\"o}lder's inequality, Sobolev's inequality, \eqref{3.3}, \eqref{3.9}, and \eqref{3.10}, we find that
\begin{align*}
|J_{1}|= & \left|-\int\rho\mathbf{u}\cdot\nabla|\mathbf{u}_t|^2dx\right| \\
\leq & 2\|\rho\|_{L^\infty}^{\frac{1}{2}}\|\mathbf{u}\|_{L^6}
\|\sqrt{\rho}\mathbf{u}_{t}\|_{L^3}\|\nabla\mathbf{u}_{t}\|_{L^2} \\
\leq & C\|\rho\|_{L^\infty}^{\frac12}\|\nabla\mathbf{u}\|_{L^2}
\|\sqrt{\rho}\mathbf{u}_{t}\|_{L^2}^{\frac12}
\|\sqrt{\rho}\mathbf{u}_{t}\|_{L^6}^{\frac12}\|\nabla\mathbf{u}_{t}\|_{L^2} \\
\leq & C\|\rho\|_{L^\infty}^{\frac34}\|\nabla\mathbf{u}\|_{L^2}
\|\sqrt{\rho}\mathbf{u}_{t}\|_{L^2}^{\frac12}
\|\nabla\mathbf{u}_{t}\|_{L^2}^{\frac32}\\
\leq & \frac{\underline{\mu}}{12}\|\nabla\mathbf{u}_{t}\|_{L^2}^2
+C\|\sqrt{\rho}\mathbf{u}_t\|_{L^2}^2;\\
|J_{2}| = & \left|-\int\rho\mathbf{u}\cdot
\nabla(\mathbf{u}\cdot\nabla\mathbf{u}\cdot\mathbf{u}_{t})dx\right| \\
\leq & \int\left(\rho|\mathbf{u}||\nabla\mathbf{u}|^2|\mathbf{u}_t|
+\rho|\mathbf{u}|^2|\nabla^2\mathbf{u}||\mathbf{u}_t|
+\rho|\mathbf{u}|^2|\nabla\mathbf{u}||\nabla\mathbf{u}_t|\right)dx \\
\leq & \|\rho\|_{L^\infty}\|\mathbf{u}\|_{L^6}
\|\nabla\mathbf{u}\|_{L^2}\|\nabla\mathbf{u}\|_{L^6}
\|\mathbf{u}_{t}\|_{L^6}
+\|\rho\|_{L^\infty}\|\mathbf{u}\|_{L^6}^2
\|\nabla^2\mathbf{u}\|_{L^2}
\|\mathbf{u}_{t}\|_{L^6} \\
\quad & +\|\rho\|_{L^\infty}\|\mathbf{u}\|_{L^6}^2
\|\nabla\mathbf{u}\|_{L^6}
\|\nabla\mathbf{u}_{t}\|_{L^2}\\
\leq & C\|\rho\|_{L^\infty}\|\nabla\mathbf{u}\|_{L^2}^2\|\mathbf{u}\|_{H^2}
\|\nabla\mathbf{u}_{t}\|_{L^2} \\
\leq & \frac{\underline{\mu}}{12}\|\nabla\mathbf{u}_{t}\|_{L^2}^2
+C\|\mathbf{u}\|_{H^2}^2;\\
|J_{3}|\leq & \|\nabla\mathbf{u}\|_{L^2}
\|\sqrt{\rho}\mathbf{u}_{t}\|_{L^4}^2
\leq C\|\nabla\mathbf{u}\|_{L^2}
\|\sqrt{\rho}\mathbf{u}_{t}\|_{L^2}^{\frac{1}{2}}
\|\sqrt{\rho}\mathbf{u}_{t}\|_{L^6}^{\frac{3}{2}} \\
\leq &C\|\rho\|_{L^\infty}^{\frac{3}{4}}\|\nabla\mathbf{u}\|_{L^2}
\|\sqrt{\rho}\mathbf{u}_{t}\|_{L^2}^{\frac{1}{2}}
\|\nabla\mathbf{u}_{t}\|_{L^2}^{\frac{3}{2}} \\
\leq & \frac{\underline{\mu}}{12}\|\nabla\mathbf{u}_{t}\|_{L^2}^2
+C\|\sqrt{\rho}\mathbf{u}_t\|_{L^2}^2;\\
|J_{4}|\leq & C\int|\mathbf{u}||\nabla\mu(\rho)||\nabla\mathbf{u}|
|\nabla\mathbf{u}_t|dx \\
\leq & C\|\mathbf{u}\|_{L^{\frac{4q}{q-2}}}\|\nabla\mu(\rho)\|_{L^q}
\|\nabla\mathbf{u}\|_{L^{\frac{4q}{q-2}}}\|\nabla\mathbf{u}_{t}\|_{L^2}
\\
\leq & C\|\nabla\mathbf{u}\|_{L^2}
\|\nabla\mathbf{u}\|_{H^1}\|\nabla\mathbf{u}_{t}\|_{L^2}
\\ \leq & \frac{\underline{\mu}}{12}\|\nabla\mathbf{u}_{t}\|_{L^2}^2
+C\|\mathbf{u}\|_{H^2}^2;\\
|J_{5}|= & \left|-\int\mathbf{b}_t\cdot\nabla\mathbf{u}_{t}\cdot\mathbf{b}dx\right|
 \\
\leq & \|\mathbf{b}_{t}\|_{L^4}\|\nabla\mathbf{u}_t\|_{L^2}
\|\mathbf{b}\|_{L^4}
\leq C\|\mathbf{b}_{t}\|_{L^2}^{\frac12}
\|\nabla\mathbf{b}_{t}\|_{L^2}^{\frac12}
\|\nabla\mathbf{u}_t\|_{L^2} \\
\leq & \frac{\underline{\mu}}{12}\|\nabla\mathbf{u}_{t}\|_{L^2}^2
+C(\delta)\|\mathbf{b}_{t}\|_{L^2}^2
+\frac{\delta}{4}\|\nabla\mathbf{b}_{t}\|_{L^2}^2;\\
|J_{6}|= & \left|-\int\mathbf{b}\cdot\nabla\mathbf{u}_{t}\cdot\mathbf{b}_t dx\right|
 \\
\leq & \|\mathbf{b}\|_{L^4}\|\nabla\mathbf{u}_t\|_{L^2}
\|\mathbf{b}_{t}\|_{L^4}
\leq C\|\nabla\mathbf{u}_t\|_{L^2}
\|\mathbf{b}_{t}\|_{L^2}^{\frac12}\|\nabla\mathbf{b}_{t}\|_{L^2}^{\frac12} \\
\leq & \frac{\underline{\mu}}{12}\|\nabla\mathbf{u}_{t}\|_{L^2}^2
+C(\delta)\|\mathbf{b}_{t}\|_{L^2}^2
+\frac{\delta}{4}\|\nabla\mathbf{b}_{t}\|_{L^2}^2.
\end{align*}
Substituting the above estimates into \eqref{3.23} and noting that
\begin{equation*}
2\int\mu(\rho)\mathfrak{D}(\mathbf{u}_t)\cdot\nabla\mathbf{u}_{t}dx
\geq \underline{\mu}\|\nabla\mathbf{u}_t\|_{L^2}^2,
\end{equation*}
we derive that
\begin{align}\label{6.2}
\frac{d}{dt}\|\sqrt{\rho}\mathbf{u}_{t}\|_{L^2}^2
+\|\nabla\mathbf{u}_{t}\|_{L^2}^2
& \leq C\|\sqrt{\rho}\mathbf{u}_t\|_{L^2}^2
+C\|\mathbf{u}\|_{H^2}^2
+\delta\|\nabla\mathbf{b}_t\|_{L^2}^2
+C\|\mathbf{b}_t\|_{L^2}^2.
\end{align}

2. It follows from \eqref{3.16} and \eqref{z3.9} that
\begin{equation}\label{3.160}
\|\mathbf{u}\|_{H^2}+\|\nabla P\|_{L^2}
\leq C\|\sqrt{\rho}\mathbf{u}_t\|_{L^2}
+C\|\nabla\mathbf{u}\|_{L^2}
+C\||\mathbf{b}||\nabla\mathbf{b}|\|_{L^2}.
\end{equation}
This along with \eqref{3.2}$_3$, \eqref{3.9}, Gagliardo-Nirenberg inequality, and Sobolev's inequality leads to
\begin{align}\label{ZX}
\|\mathbf{b}_t\|_{L^2}^2
& \leq C\|\Delta\mathbf{b}\|_{L^2}^2
+C\|\mathbf{u}\|_{L^\infty}^2\|\nabla\mathbf{b}\|_{L^2}^2
+C\|\mathbf{b}\|_{L^4}^2\|\nabla\mathbf{u}\|_{L^4}^2 \notag \\
& \leq C\|\Delta\mathbf{b}\|_{L^2}^2
+C\|\mathbf{u}\|_{L^4}\|\nabla\mathbf{u}\|_{L^4}
\|\nabla\mathbf{b}\|_{L^2}
+C\|\nabla\mathbf{b}\|_{L^2}^2
\|\nabla\mathbf{u}\|_{L^2}\|\nabla\mathbf{u}\|_{H^1}
\notag \\
& \leq C\|\Delta\mathbf{b}\|_{L^2}^2
+C\|\nabla\mathbf{u}\|_{L^2}
\|\nabla\mathbf{u}\|_{L^2}^{\frac12}\|\nabla\mathbf{u}\|_{H^1}^{\frac12}
\|\nabla\mathbf{b}\|_{L^2}
+C\|\nabla\mathbf{b}\|_{L^2}^2
\|\nabla\mathbf{u}\|_{L^2}\|\nabla\mathbf{u}\|_{H^1} \notag \\
& \leq C\|\Delta\mathbf{b}\|_{L^2}^2
+C\|\nabla\mathbf{u}\|_{H^1}^2+C\|\nabla\mathbf{u}\|_{L^2}^2
+C\|\nabla\mathbf{b}\|_{L^2}^2 \notag \\
& \leq C\big(\|\sqrt{\rho}\mathbf{u}_{t}\|_{L^2}^2
+\|\Delta\mathbf{b}\|_{L^2}^2
+\||\mathbf{b}||\nabla\mathbf{b}|\|_{L^2}^2\big)
+C\big(\|\nabla\mathbf{u}\|_{L^2}^2
+\|\nabla\mathbf{b}\|_{L^2}^2\big),
\end{align}
which combined with \eqref{6.2} and \eqref{3.160} gives
\begin{align}\label{zz3.18}
\frac{d}{dt}\|\sqrt{\rho}\mathbf{u}_{t}\|_{L^2}^2
+\|\nabla\mathbf{u}_{t}\|_{L^2}^2
& \leq C\big(\|\sqrt{\rho}\mathbf{u}_{t}\|_{L^2}^2
+\|\Delta\mathbf{b}\|_{L^2}^2
+\||\mathbf{b}||\nabla\mathbf{b}|\|_{L^2}^2\big) \notag \\
& \quad +C\big(\|\nabla\mathbf{u}\|_{L^2}^2+\|\nabla\mathbf{b}\|_{L^2}^2\big)
+\delta\|\nabla\mathbf{b}_t\|_{L^2}^2.
\end{align}

3. Differentiating \eqref{3.2}$_3$ with respect to $t$ and multiplying the resulting equations by $\mathbf{b}_t$, we obtain from integration by parts, \eqref{3.9}, and Sobolev's inequality that
\begin{equation*}
\begin{split}
\frac12\frac{d}{dt}\|\mathbf{b}_{t}\|_{L^2}^2
+\nu\|\nabla\mathbf{b}_{t}\|_{L^2}^2
& \leq C\left(\||\mathbf{u}_t||\mathbf{b}|\|_{L^2}
+\||\mathbf{u}||\mathbf{b}_t|\|_{L^2}\right)\|\nabla\mathbf{b}_{t}\|_{L^2} \\
& \leq C\left(\|\mathbf{u}_t\|_{L^4}\|\mathbf{b}\|_{L^4}
+\|\mathbf{u}\|_{L^4}\|\mathbf{b}_{t}\|_{L^4}\right)
\|\nabla\mathbf{b}_{t}\|_{L^2}
 \\
& \leq C\Big(\|\nabla\mathbf{u}_t\|_{L^2}\|\nabla\mathbf{b}\|_{L^2}
+\|\nabla\mathbf{u}\|_{L^2}
\|\mathbf{b}_{t}\|_{L^2}^{\frac12}
\|\nabla\mathbf{b}_{t}\|_{L^2}^{\frac12}\Big)
\|\nabla\mathbf{b}_{t}\|_{L^2}
\\
& \leq \frac{\nu}{2}\|\nabla\mathbf{b}_{t}\|_{L^2}^2
+C\|\nabla\mathbf{u}_t\|_{L^2}^2
+C\|\mathbf{b}_{t}\|_{L^2}^2,
\end{split}
\end{equation*}
which together with \eqref{ZX} implies that
\begin{align}\label{3.28}
\frac{d}{dt}\|\mathbf{b}_{t}\|_{L^2}^2
+\nu\|\nabla\mathbf{b}_{t}\|_{L^2}^2
& \leq C_2\|\nabla\mathbf{u}_t\|_{L^2}^2
+C\big(\|\sqrt{\rho}\mathbf{u}_{t}\|_{L^2}^2
+\|\Delta\mathbf{b}\|_{L^2}^2
+\||\mathbf{b}||\nabla\mathbf{b}|\|_{L^2}^2\big) \notag \\
& \quad +C\big(\|\nabla\mathbf{u}\|_{L^2}^2
+\|\nabla\mathbf{b}\|_{L^2}^2\big)
\end{align}
for some positive constant $C_2$.
Adding \eqref{zz3.18} multiplied by $2C_2$ to \eqref{3.28} and then choosing $\delta=\frac{\nu}{4C_2}$, we deduce that
\begin{align}\label{3.280}
& \frac{d}{dt}\Big(2C_2\|\sqrt{\rho}\mathbf{u}_{t}\|_{L^2}^2
+\|\mathbf{b}_{t}\|_{L^2}^2\Big)
+C_2\|\nabla\mathbf{u}_{t}\|_{L^2}^2
+\frac{\nu}{2}\|\nabla\mathbf{b}_{t}\|_{L^2}^2 \notag \\
& \leq C\big(\|\sqrt{\rho}\mathbf{u}_{t}\|_{L^2}^2
+\|\Delta\mathbf{b}\|_{L^2}^2
+\||\mathbf{b}||\nabla\mathbf{b}|\|_{L^2}^2\big)
+C\big(\|\nabla\mathbf{u}\|_{L^2}^2
+\|\nabla\mathbf{b}\|_{L^2}^2\big).
\end{align}

4. Multiplying \eqref{3.280} by $t^i$ ($i\in\{1,2\}$) yields
\begin{align}\label{3.18}
& \frac{d}{dt}\Big(2C_2t^i\|\sqrt{\rho}\mathbf{u}_{t}\|_{L^2}^2
+t^i\|\mathbf{b}_{t}\|_{L^2}^2\Big)
+C_2t^i\|\nabla\mathbf{u}_{t}\|_{L^2}^2
+\frac{\nu}{2}t^i\|\nabla\mathbf{b}_{t}\|_{L^2}^2 \notag \\
& \leq Ct^i\big(\|\sqrt{\rho}\mathbf{u}_{t}\|_{L^2}^2
+\|\Delta\mathbf{b}\|_{L^2}^2
+\||\mathbf{b}||\nabla\mathbf{b}|\|_{L^2}^2\big)
+Ct^i\big(\|\nabla\mathbf{u}\|_{L^2}^2
+\|\nabla\mathbf{b}\|_{L^2}^2\big) \notag \\
& \quad +C t^{i-1}\big(\|\sqrt{\rho}\mathbf{u}_{t}\|_{L^2}^2
+\|\mathbf{b}_{t}\|_{L^2}^2\big) \notag \\
& \leq Ct^i\big(\|\sqrt{\rho}\mathbf{u}_{t}\|_{L^2}^2
+\|\Delta\mathbf{b}\|_{L^2}^2
+\||\mathbf{b}||\nabla\mathbf{b}|\|_{L^2}^2\big)
+Ct^i\big(\|\nabla\mathbf{u}\|_{L^2}^2
+\|\nabla\mathbf{b}\|_{L^2}^2\big) \notag \\
& \quad +Ct^{i-1}\big(\|\sqrt{\rho}\mathbf{u}_{t}\|_{L^2}^2
+\|\Delta\mathbf{b}\|_{L^2}^2
+\||\mathbf{b}||\nabla\mathbf{b}|\|_{L^2}^2\big)
+Ct^{i-1}\big(\|\nabla\mathbf{u}\|_{L^2}^2+\|\nabla\mathbf{b}\|_{L^2}^2\big)
\end{align}
due to \eqref{ZX}. For $\sigma$ as in Lemma \ref{lem32} and any nonnegative integer $k$, we derive from \eqref{z3.3} and \eqref{z3.10} that
\begin{align}\label{7.1}
& \int_0^T t^k\big(\|\nabla\mathbf{u}\|_{L^2}^2+\|\nabla\mathbf{b}\|_{L^2}^2\big)dt
\leq \sup_{0\leq t\leq T}\big(t^ke^{-\sigma t}\big)\int_0^T e^{\sigma t}
\big(\|\nabla\mathbf{u}\|_{L^2}^2+\|\nabla\mathbf{b}\|_{L^2}^2\big)dt
\leq C,\\
& \int_0^T t^k\big(\|\sqrt{\rho}\mathbf{u}_{t}\|_{L^2}^2
+\|\Delta\mathbf{b}\|_{L^2}^2
+\||\mathbf{b}||\nabla\mathbf{b}|\|_{L^2}^2\big)dt \notag \\
& \leq \sup_{0\leq t\leq T}\big(t^ke^{-\sigma t}\big)\int_0^T e^{\sigma t}
\big(\|\sqrt{\rho}\mathbf{u}_{t}\|_{L^2}^2
+\|\Delta\mathbf{b}\|_{L^2}^2
+\||\mathbf{b}||\nabla\mathbf{b}|\|_{L^2}^2\big)dt
\leq C. \label{7.2}
\end{align}
 Integrating \eqref{3.18} over $(0,T)$ together with \eqref{7.1} and \eqref{7.2} leads to the desired \eqref{4.16}.

5. Multiplying \eqref{3.280} by $e^{\sigma t}$ together with \eqref{ZX} gives
\begin{align}\label{xxx3.18}
& \frac{d}{dt}\Big(2C_2e^{\sigma t}\|\sqrt{\rho}\mathbf{u}_{t}\|_{L^2}^2
+e^{\sigma t}\|\mathbf{b}_{t}\|_{L^2}^2\Big)
+C_2e^{\sigma t}\|\nabla\mathbf{u}_{t}\|_{L^2}^2
+\frac{\nu}{2}e^{\sigma t}\|\nabla\mathbf{b}_{t}\|_{L^2}^2 \notag \\
& \leq Ce^{\sigma t}\big(\|\sqrt{\rho}\mathbf{u}_{t}\|_{L^2}^2
+\|\Delta\mathbf{b}\|_{L^2}^2
+\||\mathbf{b}||\nabla\mathbf{b}|\|_{L^2}^2\big)
+Ce^{\sigma t}\big(\|\nabla\mathbf{u}\|_{L^2}^2
+\|\nabla\mathbf{b}\|_{L^2}^2\big) \notag \\
& \quad +C e^{\sigma t}\big(\|\sqrt{\rho}\mathbf{u}_{t}\|_{L^2}^2
+\|\mathbf{b}_{t}\|_{L^2}^2\big) \notag \\
& \leq Ce^{\sigma t}\big(\|\sqrt{\rho}\mathbf{u}_{t}\|_{L^2}^2
+\|\Delta\mathbf{b}\|_{L^2}^2
+\||\mathbf{b}||\nabla\mathbf{b}|\|_{L^2}^2\big)
+Ce^{\sigma t}\big(\|\nabla\mathbf{u}\|_{L^2}^2
+\|\nabla\mathbf{b}\|_{L^2}^2\big).
\end{align}
Integrating \eqref{xxx3.18} over $(0,T)$ together with \eqref{z3.10} and \eqref{z3.3} leads to the desired \eqref{xz4.16}.
\hfill $\Box$

\begin{lemma}\label{lem44}
Let the condition \eqref{3.1} be satisfied, then there exists a positive constant $C$ depending only on $\Omega,\underline{\mu},\bar{\mu}, \nu, q, \|\rho_0\|_{L^\infty}, \|\nabla\mathbf{u}_0\|_{L^2}^2$, and $\|\nabla\mathbf{b}_0\|_{L^2}^2$ such that
\begin{align}\label{4.29}
\int_{0}^T\|\nabla\mathbf{u}\|_{L^\infty}dt\leq C.
\end{align}
\end{lemma}
{\it Proof.}
1. Choosing $2<r<\min\{3,q\}$, we infer from Sobolev's inequality, Lemma \ref{lem23}, and \eqref{3.3} that
\begin{align}\label{10.2}
\|\nabla\mathbf{u}\|_{L^{\infty}}
& \leq C\|\mathbf{u}\|_{W^{2,r}}  \nonumber \\
& \leq C\left(\|\rho\mathbf{u}_t\|_{L^r}
+\|\rho\mathbf{u}\cdot\nabla\mathbf{u}\|_{L^r}
+\|\mathbf{b}\cdot\nabla\mathbf{b}\|_{L^r}\right)
\left(1+\|\nabla\mu(\rho)\|_{L^q}\right)^{\frac{qr}{2(q-r)}} \nonumber \\
& \leq C\|\rho\mathbf{u}_t\|_{L^3}
+C\|\mathbf{u}\|_{L^\infty}\|\nabla\mathbf{u}\|_{L^3}
+C\|\mathbf{b}\|_{L^\infty}\|\nabla\mathbf{b}\|_{L^4} \nonumber \\
& \leq C\|\rho\mathbf{u}_t\|_{L^3}+C\|\mathbf{u}\|_{H^2}^2
+C\|\mathbf{b}\|_{L^4}^{\frac12}\|\nabla\mathbf{b}\|_{L^4}^{\frac32} \nonumber \\
& \leq C\|\rho\mathbf{u}_t\|_{L^3}+C\|\mathbf{u}\|_{H^2}^2
+C\|\nabla\mathbf{b}\|_{L^2}^{\frac12}\|\nabla\mathbf{b}\|_{L^2}^{\frac34}
\|\nabla\mathbf{b}\|_{H^1}^{\frac34} \nonumber \\
& \leq C\|\rho\mathbf{u}_t\|_{L^3}+C\|\mathbf{u}\|_{H^2}^2
+C\|\nabla\mathbf{b}\|_{L^2}^2+C\|\nabla^2\mathbf{b}\|_{L^2}^2,
\end{align}
where we have used the following Gagliardo-Nirenberg inequality
\begin{align*}
\|\mathbf{b}\|_{L^\infty}\leq C\|\mathbf{b}\|_{L^4}^{\frac12}
\|\nabla\mathbf{b}\|_{L^4}^{\frac12},\
\|\nabla\mathbf{b}\|_{L^4}\leq C\|\nabla\mathbf{b}\|_{L^2}^{\frac12}
\|\nabla\mathbf{b}\|_{H^1}^{\frac12}.
\end{align*}
By H{\"o}lder's inequality, Sobolev's inequality, and \eqref{3.3}, we have
\begin{equation}\label{4.31}
\|\rho\mathbf{u}_t\|_{L^3}
\leq \|\rho\|_{L^\infty}^{\frac12}
\|\sqrt{\rho}\mathbf{u}_t\|_{L^2}^{\frac12}
\|\sqrt{\rho}\mathbf{u}_t\|_{L^6}^{\frac12}
\leq C\|\sqrt{\rho}\mathbf{u}_t\|_{L^2}^{\frac12}
\|\nabla\mathbf{u}_t\|_{L^2}^{\frac12}.
\end{equation}
As a consequence, if $T\leq1$, we obtain from \eqref{4.31} and H{\"o}lder's inequality that
\begin{align}\label{4.32}
\int_{0}^T\|\rho\mathbf{u}_t\|_{L^3}dt
& \leq C\int_{0}^T\|\sqrt{\rho}\mathbf{u}_t\|_{L^2}^{\frac12}
\|\nabla\mathbf{u}_t\|_{L^2}^{\frac12}dt \notag \\
& \leq C\Big[\int_{0}^Tt^{-\frac12}\cdot t^{-\frac13}t^{\frac13}\|\sqrt{\rho}\mathbf{u}_t\|_{L^2}^{\frac23}dt
\Big]^{\frac34}
\times\Big[\int_{0}^Tt^{\frac12}\|\nabla\mathbf{u}_t\|_{L^2}\cdot
t\|\nabla\mathbf{u}_t\|_{L^2}dt\Big]^{\frac14} \nonumber \\
& \leq C\sup_{0\leq t\leq T}\Big(t\|\sqrt{\rho}\mathbf{u}_t\|_{L^2}^2\Big)^{\frac14}
\Big(\int_{0}^Tt^{-\frac56}dt\Big)^{\frac{3}{4}}
\Big(\int_{0}^Tt\|\nabla\mathbf{u}_t\|_{L^2}^2dt\Big)^{\frac18}
\Big(\int_{0}^Tt^2\|\nabla\mathbf{u}_t\|_{L^2}^2dt\Big)^{\frac18} \nonumber \\
& \leq CT^{\frac18}\leq C.
\end{align}
If $T>1$, one deduces from \eqref{4.32} and \eqref{4.31} that
\begin{align}\label{4.33}
& \int_{0}^T\|\rho\mathbf{u}_t\|_{L^{3}}dt \notag \\
& = \int_{0}^1\|\rho\mathbf{u}_t\|_{L^{3}}dt
+\int_{1}^T\|\rho\mathbf{u}_t\|_{L^{3}}dt \nonumber \\
& \leq C+C\Big[\int_{1}^T
t^{-\frac{1}{2}}\|\sqrt{\rho}\mathbf{u}_t\|_{L^2}^{\frac{2}{3}}dt\Big]^{\frac34}
\times\Big[\int_{1}^Tt^{\frac{1}{2}}\|\nabla\mathbf{u}_t\|_{L^2}
\cdot t\|\nabla\mathbf{u}_t\|_{L^2}dt\Big]^{\frac14} \nonumber \\
& \leq C+C\Big(\sup_{1\leq t\leq T}t^2\|\sqrt{\rho}\mathbf{u}_t\|_{L^2}^2\Big)^{\frac14}
\Big(\int_{1}^Tt^{-\frac{1}{2}}\cdot t^{-\frac{2}{3}}dt\Big)^{\frac34}
\Big(\int_{1}^Tt\|\nabla\mathbf{u}_t\|_{L^2}^2dt\Big)^{\frac18}
\Big(\int_{1}^Tt^2\|\nabla\mathbf{u}_t\|_{L^2}^2dt\Big)^{\frac18} \nonumber \\
& \leq C+C\left(1-T^{-\frac{1}{6}}\right)^{\frac34}\leq C.
\end{align}
Hence, we derive the desired \eqref{4.29} from
\eqref{10.2}, \eqref{4.32}, \eqref{4.33}, \eqref{3.160}, \eqref{3.4}, and \eqref{3.9}.
\hfill $\Box$

With Lemma \ref{lem44} at hand, we immediately have the following result. The detailed proof can be found in \cite[Lemma 3.5]{Z2021} and we omit it for simplicity.
\begin{lemma}\label{lem36}
Let the condition \eqref{3.1} be satisfied, then there exists a positive constant $C$ depending only on $\Omega,\underline{\mu},\bar{\mu}, \nu, q, \|\rho_0\|_{L^\infty}, \|\nabla\mathbf{u}_0\|_{L^2}^2$, and $\|\nabla\mathbf{b}_0\|_{L^2}^2$ such that for $r\in[2,q)$,
\begin{equation}\label{4.5}
\sup_{0\leq t\leq T}\big(\|\rho\|_{W^{1,q}}+\|\rho_t\|_{L^{r}}\big)
\leq C.
\end{equation}
\end{lemma}

\begin{lemma}\label{lem35}
Let the condition \eqref{3.1} be satisfied, then there exists a positive constant $C$ depending only on $\Omega,\underline{\mu},\bar{\mu}, \nu, q, \|\rho_0\|_{L^\infty}, \|\nabla\mathbf{u}_0\|_{L^2}^2$, and $\|\nabla\mathbf{b}_0\|_{L^2}^2$ such that for $r\in(2,q)$,
\begin{align} \label{3.68}
& \sup_{0\leq t\leq T}\big[t\big(\|\mathbf{u}\|_{H^2}^2
+\|\nabla P\|_{L^2}^2+\|\mathbf{b}\|_{H^2}^2\big)\big]
+\int_0^Tt\big(\|\nabla\mathbf{u}\|_{W^{1,r}}^2
+\|\nabla P\|_{L^r}^2+\|\mathbf{b}\|_{H^3}^2\big)dt\leq C.
\end{align}
Moreover, for $\sigma$ as that in Lemma \ref{lem32} and $\zeta(T)$ as in \eqref{xz4.16}, one has
\begin{align} \label{xz3.68}
\sup_{\zeta(T)\leq t\leq T}\big[e^{\sigma t}\big(\|\mathbf{u}\|_{H^2}^2
+\|\nabla P\|_{L^2}^2+\|\mathbf{b}\|_{H^2}^2\big)\big]\leq C.
\end{align}
\end{lemma}
{\it Proof.}
1. We obtain from \eqref{3.2}$_3$, \eqref{3.9}, Sobolev's inequality, and Gagliardo-Nirenberg inequality that
\begin{align*}
\|\mathbf{b}\|_{H^2}^2
& \leq C\left(\|\mathbf{b}_t\|_{L^2}^2
+\|\mathbf{u}\cdot\nabla\mathbf{b}\|_{L^2}^2
+\|\mathbf{b}\cdot\nabla\mathbf{u}\|_{L^2}^2
+\|\mathbf{b}\|_{H^1}^2\right) \notag  \\
& \leq C\|\mathbf{b}_t\|_{L^2}^2+C\|\mathbf{u}\|_{L^4}^2
\|\nabla\mathbf{b}\|_{L^4}^2
+C\|\mathbf{b}\|_{L^\infty}^2\|\nabla\mathbf{u}\|_{L^2}^2
+C\|\nabla\mathbf{b}\|_{L^2}^2 \notag \\
& \leq C\|\mathbf{b}_t\|_{L^2}^2+C\|\nabla\mathbf{u}\|_{L^2}^2
\|\nabla\mathbf{b}\|_{L^2}\|\nabla\mathbf{b}\|_{H^1}
+C\|\mathbf{b}\|_{L^2}\|\mathbf{b}\|_{H^2}\|\nabla\mathbf{u}\|_{L^2}^2
+C\|\nabla\mathbf{b}\|_{L^2}^2 \notag \\
& \leq C\|\mathbf{b}_t\|_{L^2}^2
+\frac12\|\mathbf{b}\|_{H^2}^2
+C\|\nabla\mathbf{b}\|_{L^2}^2,
\end{align*}
which gives
\begin{align}\label{xz3.77}
\|\mathbf{b}\|_{H^2}^2
\leq C\|\mathbf{b}_t\|_{L^2}^2
+C\|\nabla\mathbf{b}\|_{L^2}^2.
\end{align}
This combined with \eqref{4.16} and \eqref{z3.10} leads to
\begin{align}\label{3.77}
\sup_{0\leq t\leq T}\big(t\|\mathbf{b}\|_{H^2}^2\big)\leq C.
\end{align}
From \eqref{3.160}, Sobolev's inequality, and \eqref{3.9}, we have
\begin{align}\label{xz3.78}
\|\mathbf{u}\|_{H^2}^2+\|\nabla P\|_{L^2}^2
& \leq C\left(\|\sqrt{\rho}\mathbf{u}_t\|_{L^2}^2+\||\mathbf{b}||\nabla\mathbf{b}|\|_{L^2}^2
+\|\nabla\mathbf{u}\|_{L^2}^2\right) \notag \\
& \leq C\left(\|\sqrt{\rho}\mathbf{u}_t\|_{L^2}^2
+\|\mathbf{b}\|_{H^2}^2\|\nabla\mathbf{b}\|_{L^2}^2
+\|\nabla\mathbf{u}\|_{L^2}^2\right) \notag \\
& \leq C\|\sqrt{\rho}\mathbf{u}_t\|_{L^2}^2
+C\|\mathbf{b}\|_{H^2}^2
+C\|\nabla\mathbf{u}\|_{L^2}^2.
\end{align}
This along with \eqref{4.16}, \eqref{3.77}, and \eqref{z3.10} yields
\begin{align}\label{3.78}
\sup_{0\leq t\leq T}\big[t\big(\|\mathbf{u}\|_{H^2}^2+\|\nabla P\|_{L^2}^2\big)\big]\leq C.
\end{align}

2. It follows from \eqref{xz3.77}, \eqref{xz4.16}, and \eqref{z3.10} that for $\zeta(T)$ as in \eqref{xz4.16},
\begin{align}\label{3.79}
\sup_{\zeta(T)\leq t\leq T}\big(e^{\sigma t}\|\mathbf{b}\|_{H^2}^2\big)\leq C,
\end{align}
which together with \eqref{xz3.78} and \eqref{z3.10} gives
\begin{align}\label{3.80}
\sup_{\zeta(T)\leq t\leq T}\big[e^{\sigma t}\big(\|\mathbf{u}\|_{H^2}^2+\|\nabla P\|_{L^2}^2\big)\big]\leq C.
\end{align}

3. For $r\in(2,q)$, we infer from Lemma \ref{lem23}, \eqref{3.3}, \eqref{3.1}, \eqref{3.9}, Sobolev's inequality, \eqref{xz3.77}, \eqref{xz3.78}, and \eqref{ZX} that
\begin{align*}
& \|\nabla\mathbf{u}\|_{W^{1,r}}^2+\|\nabla P\|_{L^r}^2 \notag \\
& \leq C\left(\|\rho\mathbf{u}_t\|_{L^r}^2
+\|\rho\mathbf{u}\cdot\nabla\mathbf{u}\|_{L^r}^2
+\|\mathbf{b}\cdot\nabla\mathbf{b}\|_{L^r}^2\right)
\left(1+\|\nabla\mu(\rho)\|_{L^q}\right)^{1+\frac{qr}{2(q-r)}}
 \nonumber \\
& \leq C\|\rho\|_{\infty}^2\|\mathbf{u}_t\|_{L^r}^2
+C\|\rho\|_{\infty}^2\|\mathbf{u}\|_{L^{\frac{qr}{q-r}}}^2
\|\nabla\mathbf{u}\|_{L^q}^2
+C\|\mathbf{b}\|_{L^{\frac{qr}{q-r}}}^2
\|\nabla\mathbf{b}\|_{L^q}^2  \nonumber \\
& \leq C\|\rho\|_{\infty}^2\|\nabla\mathbf{u}_t\|_{L^2}^2
+C\|\rho\|_{\infty}^2\|\nabla\mathbf{u}\|_{L^2}^2
\|\nabla\mathbf{u}\|_{H^1}^2
+C\|\nabla\mathbf{b}\|_{L^2}^2
\|\nabla\mathbf{b}\|_{H^1}^2 \nonumber \\
& \leq C\|\nabla\mathbf{u}_t\|_{L^2}^2
+C\big(\|\sqrt{\rho}\mathbf{u}_{t}\|_{L^2}^2
+\|\Delta\mathbf{b}\|_{L^2}^2
+\||\mathbf{b}||\nabla\mathbf{b}|\|_{L^2}^2\big)
+C\big(\|\nabla\mathbf{u}\|_{L^2}^2
+\|\nabla\mathbf{b}\|_{L^2}^2\big),
\end{align*}
which together with \eqref{4.16}, \eqref{z3.10}, and \eqref{z3.3} yields
\begin{align}\label{x2}
\int_0^Tt\big(\|\nabla\mathbf{u}\|_{W^{1,r}}^2+\|\nabla P\|_{L^r}^2\big)dt\leq C.
\end{align}
Similarly, we can show that
\begin{align}\label{x3}
\int_0^Tt\|\mathbf{b}\|_{H^3}^2dt\leq C.
\end{align}
This finishes the proof of Lemma \ref{lem35}.
\hfill $\Box$

\begin{lemma}\label{lem34}
Let the condition \eqref{3.1} be satisfied, then there exists a positive number $\varepsilon_0$ depending only on $\Omega,\underline{\mu},\bar{\mu}, \nu, q, \|\rho_0\|_{L^\infty}, \|\nabla\mathbf{u}_0\|_{L^2}^2$, and $\|\nabla\mathbf{b}_0\|_{L^2}^2$ such that
\begin{equation}\label{x}
\sup_{0\leq t\leq T}\|\nabla\mu(\rho)\|_{L^q}\leq \frac12
\end{equation}
provided that
\begin{align}\label{xy}
\|\nabla\mu(\rho_0)\|_{L^q}\leq \varepsilon_0.
\end{align}
\end{lemma}
{\it Proof.}
Taking spatial derivative $\nabla$ on the transport equation \eqref{3.10} leads to
\begin{equation}\label{6.1}
(\nabla\mu(\rho))_{t}+\mathbf{u}\cdot\nabla^2\mu(\rho)
+\nabla\mathbf{u}\cdot\nabla\mu(\rho)=\mathbf{0}.
\end{equation}
Multiplying \eqref{6.1} by $q|\nabla\mu(\rho)|^{q-2}\nabla\mu(\rho)$ and integrating the resulting equation over $\Omega$ give rise to
\begin{align*}
\frac{d}{dt}\int|\nabla\mu(\rho)|^qdx
+q\int\mathbf{u}\cdot\nabla^2\mu(\rho)\cdot
|\nabla\mu(\rho)|^{q-2}\nabla\mu(\rho)dx
 =-q\int\nabla\mathbf{u}\cdot\nabla\mu(\rho)\cdot |\nabla\mu(\rho)|^{q-2}\nabla\mu(\rho)dx.
\end{align*}
Integration by parts together with $\divv\mathbf{u}=0$ yields
\begin{align*}
q\int\mathbf{u}\cdot\nabla^2\mu(\rho)\cdot
|\nabla\mu(\rho)|^{q-2}\nabla\mu(\rho)dx
=\int\mathbf{u}\cdot\nabla(|\nabla\mu(\rho)|^q)dx
=-\int|\nabla\mu(\rho)|^q\divv\mathbf{u}dx=0.
\end{align*}
Thus, we get
\begin{align*}
\frac{d}{dt}\|\nabla\mu(\rho)\|_{L^q}^q
\leq q\int|\nabla\mathbf{u}||\nabla\mu(\rho)|^qdx
\leq q\|\nabla\mathbf{u}\|_{L^\infty}\|\nabla\mu(\rho)\|_{L^q}^q.
\end{align*}
This implies that
\begin{equation*}
\frac{d}{dt}\|\nabla\mu(\rho)\|_{L^q}
\leq\|\nabla\mathbf{u}\|_{L^\infty}\|\nabla\mu(\rho)\|_{L^q}.
\end{equation*}
which combined with Gronwall's inequality leads to
\begin{equation*}
\sup_{0\leq t\leq T}\|\nabla\mu(\rho)\|_{L^q}
\leq \|\nabla\mu(\rho_0)\|_{L^q}e^{\int_0^T\|\nabla\mathbf{u}\|_{L^\infty}dt}.
\end{equation*}
This along with \eqref{4.29} gives
\begin{align}\label{3.21}
\sup_{0\leq t\leq T}\|\nabla\mu(\rho)\|_{L^q}
&\leq C_3\|\nabla\mu(\rho_0)\|_{L^q}
\end{align}
for some constant $C_3$ depending only on $\Omega,\underline{\mu},\bar{\mu}, \nu, q, \|\rho_0\|_{L^\infty}, \|\nabla\mathbf{u}_0\|_{L^2}^2$, and $\|\nabla\mathbf{b}_0\|_{L^2}^2$. Hence, setting $\varepsilon_0=\frac{1}{2C_3}$, we obtain the desired \eqref{x} provided the condition \eqref{xy} holds true.
\hfill $\Box$

\section{Proof of Theorem \ref{thm1.1}}\label{sec4}

Suppose that the initial data $(\rho_0,\mathbf{u}_0,\mathbf{b}_0)$ satisfies \eqref{A}, according to Lemma \ref{lem20}, there exists a $T_{*}>0$ such that the problem \eqref{1.1}--\eqref{1.3} has a unique local strong solution $(\rho,\mathbf{u},\mathbf{b})$ on $\Omega\times(0,T_{*}]$. We plan to extend it to a global one. To this end, let $\varepsilon_0$ be the constant stated in Lemma \ref{lem34} and
\begin{equation}\label{5.1}
\|\nabla\mu(\rho_0)\|_{L^q}\leq \varepsilon_0.
\end{equation}

It follows from \eqref{4.5} that
\begin{equation}\label{5.2}
\rho\in C([0,T_*];W^{1,q}).
\end{equation}
Since $\mu\in C^1[0,\infty)$, we have
\begin{equation}\label{5.3}
\nabla\mu(\rho)=\mu'\nabla\rho\in C([0,T_*];L^{q}),
\end{equation}
which combined with \eqref{5.1} yields that there is a $T_1\in(0,T_*)$ such that
\begin{equation*}
\sup_{0\leq t\leq T_1}\|\nabla\mu(\rho)\|_{L^q}\leq1.
\end{equation*}
Setting
\begin{equation}\label{5.6}
T^*:=\sup\{T|(\rho,\mathbf{u},\mathbf{b})\ \text{is a strong solution on} \ \Omega\times(0,T]\},
\end{equation}
and
\begin{equation*}
T_1^*:=\sup\Big\{T\big|(\rho,\mathbf{u},\mathbf{b})\ \text{is a strong solution on} \ \Omega\times(0,T]\ \text{and}\
\sup_{0\leq t\leq T}\|\nabla\mu(\rho)\|_{L^q}\leq1\Big\}.
\end{equation*}
Then $T_1^*\geq T_1>0$. In particular, Lemmas \ref{lem33}--\ref{lem34} together with continuity arguments imply that \eqref{3.1} in
fact holds on $(0,T^*)$. Hence, we have
\begin{equation}\label{5.4}
T^*=T_1^*
\end{equation}
provided that \eqref{5.1} holds true.
Moreover, for any $0<\tau<T\leq T^*$, we infer from \eqref{3.1}, \eqref{3.9}, and \eqref{3.68} that for any $q>2$,
\begin{equation}\label{5.5}
\nabla\mathbf{u}, \nabla\mathbf{b}\in C([\tau,T];L^2\cap L^q),
\end{equation}
where one has used the standard embedding theory
$$L^\infty(\tau,T;H^1)\cap H^1(\tau,T;H^{-1})\hookrightarrow C([\tau,T];L^q)\ \text{for any}\ q\in[2,\infty).$$

Now, we claim that $T^*=\infty$. Otherwise, if $T^*<\infty$, we deduce from \eqref{5.4} that \eqref{3.1} holds at $T=T^*$.
Then it follows from \eqref{5.2}, \eqref{5.3}, and \eqref{5.5} that
$$(\rho,\mathbf{u},\mathbf{b})(x,T^*)
=\lim_{t\rightarrow T^*}(\rho,\mathbf{u},\mathbf{b})(x,t)$$
satisfies the initial conditions \eqref{A} at $t=T^*$. Thus, taking $(\rho, \mathbf{u},\mathbf{b})(x,T^*)$ as the initial data, Lemma \ref{lem20} implies that one can extend the strong solutions beyond $T^*$. This contradicts the assumption of $T^*$ in \eqref{5.6}. Moreover, the desired  exponential decay rate \eqref{1.2} follows from \eqref{xz4.16} and \eqref{xz3.68}.
This completes the proof of Theorem \ref{thm1.1}.
\hfill $\Box$\\

\section*{Acknowledgments}
The author would like to express his gratitude to the reviewers for careful reading and helpful suggestions which led to an improvement of the original manuscript.

\end{document}